\documentclass{iopart}

\usepackage{t1enc}
\usepackage{graphicx}
\usepackage[english]{babel}
\usepackage[T1]{fontenc}
\usepackage[latin1]{inputenc}
\usepackage{url}
\usepackage{graphicx}
\usepackage{graphics}
\usepackage{subfigure}
\usepackage{iopams} 
\usepackage{setstack}
\usepackage{amstext}
\usepackage{amsopn}
\usepackage{amsthm}

\newtheorem{algo}{Algorithm}

\graphicspath{{images/}}

\begin{document}

\title[Fast Sampling for L1 Problems]{Fast Markov chain Monte Carlo  sampling for sparse Bayesian inference in high-dimensional inverse problems using L1-type priors}

\author{Felix Lucka$^{1,2}$}

\address{$^1$ Institute for Computational and Applied Mathematics, University of M\"unster, Einsteinstr. 62, D-48149 M\"unster, Germany}
\address{$^2$ Institute for Biomagnetism and Biosignalanalysis, University of M\"unster, Malmedyweg 15, D-48149 M\"unster, Germany}

\ead{felix.lucka@uni-muenster.de}

\begin{abstract}
Sparsity has become a key concept for solving of high-dimensional inverse problems using variational regularization techniques. Recently, using similar sparsity-constraints in the Bayesian framework for inverse problems by encoding them in the prior distribution has attracted attention. Important questions about the relation between regularization theory and Bayesian inference still need to be addressed when using sparsity promoting inversion. A practical obstacle for these examinations is the lack of fast posterior sampling algorithms for sparse, high-dimensional Bayesian inversion:  Accessing the full range of Bayesian inference methods requires being able to draw samples from the posterior probability distribution in a fast and efficient way. This is usually done using Markov chain Monte Carlo (MCMC) sampling algorithms. In this article, we develop and examine a new implementation of a single component Gibbs MCMC sampler for sparse priors relying on L1-norms. We demonstrate that the efficiency of our Gibbs sampler  increases when the level of sparsity or the dimension of the unknowns is increased. This property is contrary to the properties of the most commonly applied Metropolis-Hastings (MH) sampling schemes: We demonstrate that the efficiency of MH schemes for L1-type priors dramatically decreases when the level of sparsity or the dimension of the unknowns is increased. Practically, Bayesian inversion for L1-type priors using MH samplers is not feasible at all. As this is commonly believed to be an intrinsic feature of MCMC sampling, the performance of our Gibbs sampler also challenges common beliefs about the applicability of sample based Bayesian inference.
\end{abstract}

\ams{65J22,62F15,65C05,65C60}
\submitto{\IP}
\maketitle


\section{Introduction} \label{sec:Intro}

\subsection{Sparse Bayesian inversion} \label{subsec:SpBayInv}
Solving high-dimensional inverse problems using \emph{sparsity} constraints as \emph{a priori information} has led to enormous advances in various application areas. \emph{Total variation} (\emph{TV}) deblurring \cite{RuOsFa92,BuOs12} uses sparsity constraints on the gradient of the unknown quantity and is successfully used in many imaging applications. By the notion of \emph{compressed sensing} \cite{Do06,ElKu12}, a number of techniques are summarized, which rely on the idea that high quality reconstructions can be obtained from a small amount of data, if a sparse basis for the unknowns is a priori known. Traditionally, sparsity constraints are formulated in the framework of \emph{variational regularization} when introduced as a priori information to inverse problems. One popular approach is to use regularization functionals incorporating L1 norms. However, this type of sparsity constraints can also be formulated and examined in the framework of Bayesian statistics \cite{KoLaNiSi12}. Addressing high dimensional, ill-posed inverse problems as problems of Bayesian inference has gained growing attention over the years \cite{TaVa82,GeGe84,KaSo05,Id08}. It allows an easy formulation of \emph{a priori information} on the solution via \emph{a priori probability distributions} (\emph{prior}). Furthermore, a specific inference strategy that is called the \emph{maximum a posteriori} estimate (\emph{MAP}) corresponds to variational regularization (the prior corresponds to the regularization functional). In general, the Bayesian solution to an inverse problem is given by the \emph{a posteriori probability distribution} (\emph{posterior}) over the parameter space (the MAP estimate is the point maximizing this distribution). The analysis of sparse Bayesian inversion is far less elaborate up to now, and a number of exciting questions still remain to be addressed. In particular, sparse inversion is an interesting topic to study the relation between regularization theory and Bayesian inference. This relation is well understood for regularization using L2 norms, which corresponds to Bayesian inference with Gaussian priors, see, e.g., \cite{KaSo05}. While the differences in these scenarios are subtle, they become way more pronounced in the context of sparse inversion using \emph{L1-type priors}, i.e., priors, which rely on L1 norms of the unknowns \cite{LaSi04,LaSaSi09}. A central tool to study these differences is the examination of the posterior by \emph{Monte Carlo sampling} methods. The standard sampling techniques were designed for Bayesian inference in low-dimensional, well-posed problems and often fail when used in scenarios arising from typical inverse problems. For these reasons, a number of specific sampling techniques for ill-posed, high-dimensional problems have already been developed \cite{HiLeHo03,HaLaLeSaTa04,HaLaMiSa06,MoFoSv08,PaFo12}. However, these techniques mainly address Gaussian priors. For the sparsity-promoting  L1-type priors they may fail dramatically. This observation is in line with the fact that efficient optimization techniques for the corresponding variational regularization schemes are still a vital field of research as well \cite{GoOs09,BuMoBeOs11}.

\subsection{Contributions and Structure} \label{subsec:ConStr}

Our work was motivated by questions that arise when Bayesian inference using general L1-type priors is applied to edge preserving image reconstruction, similar to the scenarios discussed in \cite{LaSi04,LaSaSi09,KoLaNiSi12}. A major problem we and others faced was that the conventional \emph{Markov chain Monte Carlo} (\emph{MCMC}) tools for sampling the posterior distribution fail in such situations, rendering many examinations infeasible. Therefore, we develop and examine new and more efficient implementations of sampling algorithms for these situations first. In this paper, we present a fast implementation of a Gibbs sampling algorithm and study its performance in two typical inverse problems scenarios. Thereby, we provide a solid basis for addressing more sophisticated questions in sparse Bayesian inversion in the future.\\ 
In Section \ref{sec:Methods}, we describe the setting and methods used. Detailed numerical examinations of all MCMC algorithms for two test scenarios are presented in Section \ref{sec:Results}. In Section \ref{sec:Discussion}, the results are discussed and we point to future directions of development. Additionally, implementation  details and code for the new sampling algorithms is provided in \ref{sec:App_Code}.


\section{Methods} \label{sec:Methods}
In this section, we will first introduce the general setting for our examinations (Section \ref{subsec:GenSetBayFor}) and the basics of Bayesian inference. Then we review the basic principles of MCMC-based posterior inference and present the most popular MCMC sampling schemes, the \emph{Metropolis-Hastings} algorithm and the \emph{Gibbs sampling} algorithm (Section \ref{subsec:PosInfMCMSam}). The intention of these first two sections is to make the article more accessible for readers which manly used variational regularization techniques so far and have little experience with Bayesian techniques. The more advanced reader can skip these sections. Section \ref{subsec:ImpGibbsL1} contains the main contributions of this article, i.e., the development of a \emph{Gibbs sampling} scheme for L1-type priors which relies on a robust numerical implementation of an exact, explicit sampling from the conditional single component posterior. In the last section (Section \ref{subsubsec:MCMCConDia}), we will explain the methods used for the evaluation of the sampling performance in the computational studies. The experienced reader may, again, skip this section.

\subsection{General Setting and Bayesian Formulation} \label{subsec:GenSetBayFor}

In general, we consider the inverse problem of solving a continuous, linear, ill-posed operator equation. Here, we start from the following discrete model chosen for obtaining a computational solution (the \emph{computational model}):
 \begin{equation}
  m = A \, u + \varepsilon, \label{eq:FwdEq1}
 \end{equation}
where $m \in \mathbb{R}^k$ represents the given measurement data, $u \in \mathbb{R}^n$ represents the unknowns derived from a discretization of the computational domain, $ A \in \mathbb{R}^{k \times n}$ is the discretization of the continuous forward operator with respect to the domains of $u$ and $m$ and $\varepsilon  \in \mathbb{R}^k$ is an additive, stochastic noise term. Accounting for the stochastic nature of the noise term renders \eref{eq:FwdEq1} into a relation between the $k$-dim random variables $M$ and $\mathcal{E}$ (the \emph{likelihood model}):
\begin{equation}
 M =   A \,  u + \mathcal{E} \label{eq:LikeMod}
\end{equation}
See  \cite{Ho06,BiHoMuRu07} for details on the implications of this step. For simplicity, we assume $\mathcal{E} \sim \mathcal{N}(0,\sigma^2 I_k),  \sigma > 0$ here, where $I_k$ is the $k$-dim identity matrix (the extension to general Gaussian noise is straight forward). Now, the conditional probability density of $M$ given $u$ is determined by \eref{eq:LikeMod} and is, thus, called the \emph{likelihood} density:
\begin{equation}
  p_{li}(m|u) = \left( \frac{1}{2 \pi \sigma^2} \right)^{\case{k}{2}} \exp \left( -\frac{1}{2 \, \sigma^2} \| m - A \, u \|^{2}_{2}\right) \label{Likelihood}
\end{equation}
Due to the ill-posedness of \eref{eq:FwdEq1}, inference about $u$ given $M$ on the basis of \eref{Likelihood} is not feasible with standard statistical inference strategies. \emph{Bayesian inference strategies} rely on considering $u$ as  a random variable itself ($U$ in our notation) and on encoding a priori information about $U$ in its density, $p_{pr}(u)$, which is therefore called the \emph{prior}. Then, the model can be inverted using Bayes' rule:
\begin{equation}
  p_{post}(u|m) = \frac{p_{li}(m|u)p_{pr}(u)}{p(m)} \label{eq:BayesRule}
\end{equation}
The conditional density of $U$ given $M$ is called the \emph{posterior}. In Bayesian inference, this density is the complete solution to the inverse problem. The term $p(m)$ is called the \emph{model-evidence} and for our aims, it is just a normalizing constant, which is of no further importance. There are several ways to exploit the information about $U$ contained in the posterior. The most popular one, called the \emph{maximum a posteriori} estimate (\emph{MAP}), is to infer a point estimate for $U$ by searching for the highest \emph{mode} of the posterior. Another way to obtain a point estimate, called the \emph{conditional mean} estimate (\emph{CM}), is to compute the mean/expected value of the posterior:
\begin{eqnarray}
\hat{u}_{\rm MAP} &:= \underset{{u \in \mathbb{R}^n}}{{\rm argmax}} \left\lbrace   \; p_{post}(u|m) \right\rbrace \\
\hat{u}_{\rm CM} &:= \mathbb{E} \left[ u|m \right] = \int u \; p_{post}(u|m) \; \rmd u \label{eq:CMDef}
\end{eqnarray}
Practically, computing the MAP estimate is a high-dimensional \emph{optimization} problem, whereas computing the CM estimate is a high-dimensional \emph{integration} problem. Apart from point estimates, computing \emph{confidence intervals}, \emph{conditional covariance} or \emph{histogram} estimates are other applications of posterior-based inference. See \cite{KaSo05} for an overview and, e.g., \cite{HaLaLeSaTa04,HaLaMiSa06,CaSo07} for the applications to remote sensing, algae population dynamics and image deblurring. \\
This far, we did not specify the concrete form of the prior $p_{pr}(u)$, which is actually the most important step within the Bayesian formalism. A common choice linking Bayesian inference with variational regularization is given by \emph{Gibbs distributions}:
\begin{equation}
  p_{pr}(u) \propto \exp \left(- \lambda \mathcal{J}(u) \right)  \label{eq:GibbsPrior}
\end{equation}
Here, $\mathcal{J}(u)$ is an energy functional penalizing unwanted features of $u$, and $\lambda > 0$ is a scaling parameter that is called the \emph{regularization parameter}. Now, after suppressing terms not dependent on $u$, the MAP estimate is given by
\begin{eqnarray}
\hat{u}_{\rm MAP} &= \underset{{u \in \mathbb{R}^n}}{{\rm argmax}} \; \left\lbrace  \exp \left( -\frac{1}{2 \, \sigma^2} \| m - A \, u \|^{2}_{2} -  \lambda \mathcal{J}(u)\right) \right\rbrace  \nonumber \\
&=\underset{{u \in \mathbb{R}^n}}{{\rm argmin}} \left\lbrace  \|  m - A \, u  \|^{2}_{2} + (2 \, \sigma^2 \lambda ) \; \mathcal{J}(u) \right\rbrace \label{eq:Tikh}
\end{eqnarray}
This is a \emph{Tikhonov-type} regularization of equation \eref{eq:FwdEq1} \cite{EnHaNe96}.\\ 
In this article, we only consider Gibbs priors with a L1 norm type energy functional:
\begin{equation}
  \mathcal{J}(u)  = | D \, u|,
\end{equation}
where $D \in  \mathbb{R}^{l \times n}$, and $|\cdot|$ denotes the L1 norm in $\mathbb{R}^l$. Although such priors may seem like a generic extension of the one dimensional \emph{Laplace distribution} to a multivariate setting, we note here that multivariate generalizations of Laplace distributions are commonly defined in a different way (see, e.g. \cite{ElAtTe06}). Concrete examples of L1-type priors will be given in Section \ref{sec:Results}. For the methods presented here, we require $l \leqslant n$, ${\rm rank}(D) = l$ and $\ker (D) \cap \ker (A) = {0}$. In forthcoming work, we will extend the sampler proposed in Section \ref{subsec:ImpGibbsL1} to more general settings. \\

\textit{Remark:} For the sake of an intuitive presentation of the Bayesian formulation of inverse problems, we started from a deterministic setting. For a more detailed and rigorous description on how to derive a discrete, computational model of the continuous inverse problem starting in the Bayesian framework, we refer to \cite{LaSaSi09}.

\subsection{Posterior Inference using MCMC Sampling} \label{subsec:PosInfMCMSam}

\paragraph{General Principles:} \label{para:GenPri}

In typical inverse problems scenarios, the dimension $n$ of the unknowns $u$ is very large (in our computational examples, we will study a scenario where the limit $n \rightarrow \infty$ is of central interest). Therefore, the integration to compute the CM estimate \eref{eq:CMDef} is intractable by means of traditional quadratures. Interval, conditional covariance and histogram estimates and even more sophisticated topics in Bayesian inference like \emph{marginalization}, \emph{model selection} or \emph{experiment design} \cite{To11} also rely on integration tasks and can, thus, not be computed by such an approach as well. Integration by \emph{Monte Carlo} methods can avoid these difficulties. A sequence of points $u_i,\, i=1,\ldots,K$ is constructed, which is distributed like the posterior (the construction schemes are called \emph{sampler} or \emph{sampling methods}). If they were drawn independently, the \emph{law of large numbers} would guarantee that
 \begin{equation}
   \frac{1}{K} \sum_{i=1}^K f(u_i) \, \overset{K \rightarrow \infty}{\longrightarrow} \, \mathbb{E} \left[ f(u)|m \right] = \int_{\mathbb{R}^n}f (u) \, p_{post}(u|m) \: \rmd u \label{eq:LLN}
  \end{equation}
for any measurable $f$ almost surely and in L1 with rate $O (K^{-1/2})$. This means that the empirical mean of the sequence $f(u_i)$, $i=1,\ldots,K$ converges to the expected value of $f(u)$ w.r.t the posterior \cite{Kl08}. A difficulty in our setting is that the posterior is not given in a form that allows for drawing independent samples. It is only known up to a normalizing constant (the model-evidence) and does not belong to a class of distributions for which independent sampling schemes are known. However, by the \emph{strong ergodic theorem}, the above convergence (and its rate) still holds if the sequence is dependent, but originates from an \emph{ergodic Markov chain} that has $p_{post}(u|m)$ as its \emph{equilibrium distribution} \cite{Kl08}. Techniques to construct such chains are called \emph{Markov chain Monte Carlo} (\emph{MCMC}) methods. A huge number of different MCMC methods have been proposed. However, no method is known, which exhibits  a good performance for all types of distributions. For a comprehensive overview, we refer to \cite{Li08}, for the application to inverse problems, see \cite{KaSo05}. Most MCMC methods rely on one of two basic sampling schemes, which we will introduce in the next sections. Instead of comparing all possible and sophisticated variants of these schemes in our studies, we will use a small number of simple variants and focus on the differences between the two basic schemes for L1-type priors.

\paragraph{Metropolis-Hastings Sampling:} \label{para:MetHasSch}

For the ease of presentation, we denote the target probability density we want to sample by $p(x)$, $x \in \mathbb{R}^n$. The \emph{Metropolis-Hastings} (\emph{MH}) algorithm \cite{MeRoRoTeTe53,Ha70} is a very simple rule to generate a Markov chain:
\begin{algo}{\textbf{(Metropolis-Hastings Sampling)}}
Let $q(x,y): \mathbb{R}^n \times \mathbb{R}^n \rightarrow \mathbb{R}_+$ be a function satisfying $\int q(x,y) \rmd y = 1$ for all $x \in \mathbb{R}^n$ (proposal distribution) and $x_0 \in \mathbb{R}^n$ an initial state. Define burn-in size $K_0$ and sample size $K$.\\
For $i$ $=$ $1$,$\ldots$,$K_0+K$ do:
 \begin{enumerate}
  \item[1] Draw $y$ from the proposal distribution $q(x_{i-1},y)$.
  \item[2]  Compute the acceptance ratio 
  \begin{equation*}
   	r(x_{i-1},y) = \min \left(1,  \frac{p(y)\, q(y,x_{i-1})}{p(x_{i-1}) \, q(x_{i-1},y)} \right).\label{eq:AccRatio}
  \end{equation*}
  \item[3]  Draw $\theta \in [0,1]$ from a uniform probability density.
  \item[4]  If $r \geqslant \theta$, set $x_i = y$, else  set $x_i = x_{i-1}$.
 \end{enumerate}
 Return $x_{K_0+1},\ldots,x_{K}$.
\end{algo}
Note that the restrictions on $p(x)$ for this scheme are minimal: We only have to know $p(x)$ up to a scaling factor, as only ratios of probabilities are used, and we only need to be able to evaluate $p(x)$ for any given $x$. Each sampling step requires one such evaluation (in inverse problems, the computational demanding part of this evaluation is usually applying the forward mapping $A$). The numerical implementation of the raw MH scheme is trivial. MH can, thus, be considered as a ``black-box sampler'', which explains its success in many different application areas \cite{Li08}.\\
However, while the scheme works for all kinds  of proposal distributions in theory, its application is only feasible if $q(x,y)$ leads to a chain that moves ``fast'' in the sampling space with respect to computational speed. This way, the important regions of the sampling space are explored reasonably fast and consecutive samples are as uncorrelated as possible (which improves the convergence in \eref{eq:LLN}). These requirements are hard to fulfill in practice. We refer to the discussions in \cite{KaSo05,Li08}.
Usually, one ends up in a well-known dilemma of tuning different opposing parameters by manual inspection of different chain characteristics (see Section \ref{subsubsec:MCMCConDia}). Additionally, different applications usually require to develop and implement specific proposal distributions. As a consequence, a huge number of different MH-based schemes exist \cite{Li08}. Especially for inverse problems, sophisticated algorithms that include automatic tuning procedures for the sampling parameters have been developed \cite{HaLaLeSaTa04,HaLaMiSa06,ChFo10}. However, as mentioned earlier, a detailed comparison of their performance in sparse Bayesian inversion is not the topic of this publication. We rather want to compare the basic variants of MH algorithms and their performance to basic variants of Gibbs sampling algorithms. Therefore, we will use three proposal distributions that are commonly applied in practice because of their simplicity. They all belong to the class of \emph{symmetric random-walk} Metropolis schemes \cite{Li08}:
\begin{equation}
 y = x + \vartheta, \qquad  \mathbb{E}(\vartheta) = 0, \qquad p_\vartheta (\omega) \varpropto g( \| \omega \|_2) \, \forall \, \omega \in \mathbb{R}^n, 
\end{equation}
for a suitable, non-negative function $g$. This means that a new proposal $y$ is generated by perturbing the current state $x$ in a random, unbiased, symmetric way. Thus, $q(x,y) \varpropto g(\| x-y \|)$, and $q$ vanishes from the acceptance ratio \eref{eq:AccRatio}. The three choices of $\vartheta$ we will use are:
\begin{description}
 \item[1) MH-Iso:] All components of $x$ are updated: $\vartheta_i \sim \mathcal{N}(0,\kappa^2), \, \forall \, i$.
 \item[2) MH-Ncom:] $1<n^* < n$ components $i_1,\ldots,i_{n^*}$ of $x$ are randomly chosen and are updated while the other components remain unchanged: $\vartheta_i \sim  \mathcal{N}(0,\kappa^2)$. if $i \in \{i_1,\ldots,i_{n^*}\}$, else $\vartheta_i = 0$.
  \item[3) MH-Si:] One component $i_*$ of $x$ is randomly chosen and updated while all other components remain unchanged: $\vartheta_{i_*} \sim  \mathcal{N}(0,\kappa^2)$, $\vartheta_{[-i_*]} = 0$.
 \end{description}
 Here, $\vartheta_{[-i]}$ denotes all components of $\vartheta$ except the $i^{th}$ one. The concrete choice of $n^*$ and $\kappa$ will be explained in Section \ref{sec:Results}.

\paragraph{Gibbs Sampling:} \label{para:GibbsSam}

In certain scenarios, direct sampling of a $n$-dim multivariate distribution is not possible or computationally too expensive, but direct sampling from conditioned (thus, lower dimensional) versions of that distribution is feasible. In such a situation, \emph{Gibbs sampling} can be applied. By successive sampling from the lower dimensional conditional distributions while changing the coordinates, which are fixed in each step, a Markov chain is generated  \cite{GeGe84,GeSm90}. The most basic scheme is given by:
\begin{algo}{\textbf{(Single Component Gibbs Sampling)}} \label{algo:Gibbs}
Let $x_0 \in \mathbb{R}^n$ an initial state. Define burn-in size $K_0$ and sample size $K$\\
For $i$ $=$ $1$,$\ldots$,$K_0+K$ do:
\begin{itemize}
 \item Set $x_i := x_{i-1}$.
 \item[] For $j$ $=$ $1$,$\ldots$,$n$ do:
 \begin{enumerate}
  \item[1] Set $s = j$ (systematic scan) or draw $s$ randomly from $\left\lbrace 1,\ldots,n \right\rbrace $ (random scan).
  \item[2] Draw $(x_i)_{s}$ from the conditional, 1-dim density $p(\, \cdot \, |(x_i)_{[-s]})  $.
 \end{enumerate}
 \end{itemize}
 Return $x_{K_0+1},\ldots,x_{K}$.
\end{algo}
We will abbreviate the systematic scan version of the above sampler as \textbf{SysGibbs} and the random version (which requires the extra computational effort of picking a random coordinate) as \textbf{RnGibbs}.\\
The basic Gibbs sampling scheme can be very slow if the correlations between the single components $x_i$ are strong. This occurs naturally in typical under-determined inverse problems. In this case, the conditional distributions differ considerably from the corresponding marginal ones. As a consequence, the chain moves very randomly, exploring the search space very slowly. To address this problem, \emph{overrelaxed} variants of Gibbs sampling have been proposed. The idea behind them are similar to those used in overrelaxation techniques for the iterative solution of systems of linear equations \cite{Sa03}. The specific form of overrelaxation that we will apply and examine was proposed in \cite{Ne95}, and relies on \emph{order statistics}. Step 2 in Algorithm \ref{algo:Gibbs} is replaced by:
\begin{algo}{\textbf{(Ordered Overrelaxation)}} \label{algo:OOR}
  \begin{enumerate}
  \item[2.1] Draw $N_O$ random values from the conditional, 1-dim density $p(\, \cdot \, |(x_i)_{[-s]})$, where $N_O \in \mathbb{N}$ is odd.
  \item[2.2] Arrange these $N_O$ values plus the old value $(x_i)_s$ in non-decreasing order, labeling them as follows:
  \begin{equation}
   (x_i)_s^{(0)} \leqslant (x_i)_s^{(1)} \leqslant \cdots \leqslant (x_i)_s^{(t)} = (x_i)_s \leqslant \cdots \leqslant (x_i)_s^{(N_O)}
  \end{equation}
  \item[2.3] Replace $(x_i)_s$ by $(x_i)_s^{(N_O - t)}$.
 \end{enumerate}
\end{algo}
The value of $N_O$ functions like an overrelaxation parameter. The larger the value of $N_O$, the larger the effect of overrelaxation and more randomness of the sampling process is suppressed. For symmetric densities, the current value of the component is mirrored at the mean and the whole chain moves on an iso-probability level of the density in the limit of $N_O \rightarrow \infty$. We will discuss more details of ordered overrelaxation in Section \ref{subsubsec:NormVsOOR}. We will denote the overrelaxed versions of SysGibbs and RnGibbs by appending "O$N_O$", e.g., "SysGibbsO7" denotes the systematic scan Gibbs Sampler with ordered overrelaxation using $N_O = 7$. While the basic scheme for ordered overrelaxation requires $N_O$ times more computation time compared to Algorithm \ref{algo:OOR}, an efficient implementation is given in \cite{Ne95} that renders the computation time nearly independent of $N_O$. We will present this form after the next paragraph.

\subsection{Implementation of Gibbs Sampling for L1-type Priors} \label{subsec:ImpGibbsL1}

In this section, the main contributions of this article are presented. The general Gibbs sampling schemes (Algorithms \ref{algo:Gibbs} and \ref{algo:OOR}) need to be implemented in an efficient way. For this, we will first derive a way to compute a simple representation of the conditional single component density and then explain how to implement an exact, explicit and numerically robust sampler for it.

\paragraph{Conditional Densities for L1-type Priors:} \label{para:SinCompDensL1}
We will now derive the single component conditional densities required by Algorithms \ref{algo:Gibbs} and \ref{algo:OOR} for our setting (cf. Section \ref{subsec:GenSetBayFor}). Because ${\rm rank}(D) = l$ (cf. Section \ref{subsec:GenSetBayFor})., we can find $v_1,\ldots,v_l \in \mathbb{R}^n$ such that $D \, v_i = e_i $, (where $e_i$ denotes the $i^{th}$ unit vector in $\mathbb{R}^l$) and $v_{l+1},\ldots,v_n \in \mathbb{R}^n$ such that $v_1,\ldots,v_n$ form a basis of $\mathbb{R}^n$. Then, we have
\begin{equation*}
 \fl \quad u = \sum_{i = 1}^{n} \xi_i v_i \;, \quad  \quad D \, u = (\xi_1,\ldots,\xi_l)^t \quad \quad {\rm and} \quad \quad | D \, u | = \sum_{i = 1}^{l} | \xi_i |.
\end{equation*}
With $V := [v_1,\ldots,v_n]$, we can transform the posterior to:
\begin{eqnarray}
 \fl \quad \exp \left(  - \frac{1}{2 \, \sigma^2} \| m - A \, u  \|^{2}_{2} \,-\, \lambda | D \, u | \right)  & = & \exp \left( - \frac{1}{2 \, \sigma^2} \|  m - A \, V  \xi \|^{2}_{2} \,-\, \lambda \sum_{i = 1}^{l} | \xi_i | \right)  \nonumber \\
& := & \exp \left(  - \|  \bar{m} - \Psi \, \xi \|^{2}_{2} \,-\, \lambda \sum_{i = 1}^{l} | \xi_i | \right) ,
\end{eqnarray}
where $\bar{m} := m/(\sqrt{2} \, \sigma)$ and $\Psi := (A \, V) /(\sqrt{2} \, \sigma)$. Because the transformations are linear, no specific attention to the correct transformation of probability densities has to be paid. Now let $\psi_i$ be the $i^{th}$ column of $\Psi$, $\Psi_{[-i]}$ be $\Psi$ without the $i^{th}$ column and $\xi_{[-i]}$ be $\xi$ without the $i^{th}$ entry. Then
\begin{equation}
  \bar{m} -\Psi \, \xi =  (\bar{m} - \Psi_{[-i]} \, \xi_{[-i]}) - \psi_i \, \xi_i :=  \varphi_{[-i]} -\psi_i \, \xi_i.
\end{equation}
Consider the conditional posterior of $\xi_i$ given $\bar{m}$ and $\xi_{[-i]}$:
\begin{eqnarray}
\fl \qquad p(\xi_i|\bar{m},\xi_{[-i]}) &\propto \exp\left(- \|  \varphi_{[-i]} -\psi_i \, \xi_i  \|^{2}_{2} -  \lambda  | \xi_i |  \cdot \mathbf{1}_{\{ i \leqslant l \}}\right)  \nonumber \\
			     &= \exp\left(- \left\langle   \varphi_{[-i]} -\psi_i \, \xi_i,  \varphi_{[-i]} -\psi_i \, \xi_i \right\rangle  -  \lambda  | \xi_i | \cdot \mathbf{1}_{\{ i \leqslant l \}} \right)  \nonumber \\
			     &\propto \exp\left(- \|  \psi_i \|^{2}_{2} \, \xi_i^2 + 2 \psi_i^t \varphi_{[-i]} \, \xi_i  -  \lambda | \xi_i |  \cdot \mathbf{1}_{\{ i \leqslant l \}} \right)
\end{eqnarray}
To ease the following presentation, we define:
\begin{equation}
 \fl  x := \xi_i; \; a := \|  \psi_i \|^{2}_{2}; \: b := 2 \psi_i^t \varphi_{[-i]} = 2 \left[  \psi_i^t \, \bar{m} - \left( \psi_i^t \, \Psi_{[-i]}\right)  \, \xi_{[-i]}\right]; \: c := \lambda \cdot \mathbf{1}_{\{ i \leqslant l \}} \label{eq:ABC}
\end{equation}
Thus, the problem of sampling from the single component conditional densities is reduced to sampling from the 1-dim density 
\begin{equation}
 p(x) \propto \exp (- a \, x^2 + b \, x - c \, |x|) \label{StdDens},
\end{equation}
once $a$, $b$ and $c$ have been computed by \eref{eq:ABC}. In the next paragraph, we will describe how to use the \emph{inverse cumulative distribution method} \cite{KaSo05} to sample from \eref{StdDens}. Concerning the practical implementation of computing $a$, $b$ and $c$ in a fast way, note that $a$ and $c$ can be precomputed and only $b$ depends on the current state of the chain $\xi$ through the term $(\psi_i^t \, \Psi_{[-i]}) \, \xi_{[-i]}$. The most efficient way to compute this term strongly depends on the form of $A$ and $V$, on the problem size $n$ and on the hardware available. If enough working memory is available to store the $n \times n$ matrix $\Phi := \Psi^t  \Psi$, the most efficient way is to compute 
\begin{equation}
 (\psi_i^t \, \Psi_{[-i]}) \, \xi_{[-i]} = \xi^t \Phi_{(\cdot,i)} - \xi_i  \|  \psi_i \|^{2}_{2},
\end{equation}
because the most extensive operation is a scalar product of dimension $n$. In the scenario examined in \cite{LaSaSi09}, $A$ is a symmetric convolution operator and $V$ an inverse wavelet transform, i.e., $v_j$ are the wavelets. For large $k$ and $n$ (as encountered, e.g., in 2D or 3D imaging applications), it is infeasible to compute and store the matrix form of $A$, $V$ or $\Psi$. Then, it is advantageous to use
\begin{eqnarray}
 \fl  (\psi_i^t \, \Psi_{[-i]}) \, \xi_{[-i]} = \psi_i^t  \left( \Psi \xi \right) - \xi_i  \|  \psi_i \|^{2}_{2} = \frac{1}{2 \sigma^2} \left( A \cdot V e_i   \right)^t  \left( A V \xi   \right)    - \xi_i  \|  \psi_i \|^{2}_{2} \nonumber \\
 \fl  \qquad = \frac{1}{2 \sigma^2} \left[ A \cdot \left(V e_i \right)  \right]^t  \left[ A \left( V \xi \right)  \right]    - \xi_i  \|  \psi_i \|^{2}_{2} =  \frac{1}{2 \sigma^2} v_i^t  \left[ \left(A^t A \right) \cdot \left(V \xi \right)   \right]  - \xi_i  \|  \psi_i \|^{2}_{2}. \label{eq:2DimplB}
\end{eqnarray}
Here, $\left(V \xi \right) $ can be realized using the \emph{fast wavelet transform}, while the double convolution by $ \left(A^t A \right)$ can be substituted by a single convolution with a different kernel and realized by the \emph{fast Fourier transform}.

\paragraph{Explicit 1D Sampling:} \label{para:Exp1DSam}

Sampling from continuous 1-dim distributions by the inverse cumulative distribution method follows a simple rule: Let $F(y) := \int_{-\infty}^y p(x) \rmd x$ be the \emph{cumulative distribution function} (\emph{cdf}) and $r$ be a random number uniformly drawn from $[0,1]$. Then, $y = F^{-1}(r)$ is distributed like $p(x)$ (see \cite{Li08}). We can also use this concept to provide an equivalent implementation of Algorithm \ref{algo:OOR}. For for a given $N_O$:
\begin{algo}{\textbf{(CDF Implementation of Ordered Overrelaxation)}} \label{algo:ImplOOR}
 \begin{enumerate}
  \item[2.1] Compute $r = F\left[ (x_i)_s\right] $, which lies in $[0,1]$.
  \item[2.2] Let $r'$ be the random ordered overrelaxation of $r$ w.r.t to the uniform distribution on $[0,1]$ and $N_O$ (computed with Algorithm \ref{algo:OOR}).
  \item[2.3] Replace $(x_i)_s$ by $F^{-1}(r')$.
 \end{enumerate}
\end{algo}
For more details, we refer to \cite{Ne95}. Turning these rules into efficient sampling schemes requires a fast and stable way to invert $F$, which is defined by an integral. Using numerical integration for this purpose often fails to render fast and robust sampling algorithms. In this paragraph, we will present a scheme that relies on the \emph{inverse complementary error function} (erfcinv) for which efficient and stable implementations are known. \\
First, we compute the normalization factor for $p(x) \propto \exp (- a \, x^2 + b \, x - c \, |x|)$. Splitting the integral from $-\infty$ to $\infty$ into two parts (from $-\infty$ to $0$ and the rest) yields subproblems that can be treated like the normalization of the normal distribution (completing the square and a linear integral transformation). This leads to:
\begin{eqnarray}
 \mathcal{N} &:= \int_{-\infty}^{\infty} \exp (- a \, x^2 + b \, x - c \, |x|) \rmd x  \nonumber \\
&= \frac{1}{2} \sqrt{\frac{\pi}{a}} \left[  e^{\frac{(b+c)^2}{4a}} {\rm erfc}\left( \frac{b+c}{2\sqrt{a}}\right)  + e^{\frac{(c-b)^2}{4a}} {\rm erfc}\left( \frac{c-b}{2\sqrt{a}} \right) \right]  \nonumber \\
&:= \chi \left[ \tilde{e}_+ \, {\rm erfc}\left( \alpha_+ \right)  + \tilde{e}_- \, {\rm erfc}\left( \alpha_- \right) \right],
\end{eqnarray}
where ${\rm erfc}(y) := \frac{2}{\sqrt{\pi}} \int_{y}^{\infty} e^{-t^2} \, \rmd t$ denotes the complementary error function. The cdf is given by:
\begin{eqnarray}
 \fl &{\rm cdf}(y) := \frac{1}{\mathcal{N}} \int_{-\infty}^y \exp (- a \, x^2 + b \, x - c \, |x|) \rmd x \nonumber \\
 \fl &\qquad \;\; = \frac{\chi}{\mathcal{N}}  \cdot \cases{\tilde{e}_+ {\rm erfc} \left( - \sqrt{a} \, y + \alpha_+ \right) &, if y < 0,\\
       \tilde{e}_+ {\rm erfc} \left( \alpha_+ \right) + \tilde{e}_- \left[ {\rm erfc} \left( \alpha_- \right) - {\rm erfc}\left( \sqrt{a} \, y + \alpha_-  \right)  \right] &, if y > 0.\\} \label{eq:cdf}
\end{eqnarray}
Inverting this cdf for a given $r \in [0,1]$ is simple. To find $y = {\rm cdfinv} (r)$ we first check if $y < 0$  by using the cdf for this domain. Let
\begin{eqnarray}
\fl z := {\rm erfcinv} \left( \frac{r \, \mathcal{N}}{\chi \, \tilde{e}_+}  \right) &= {\rm erfcinv} \left\lbrace  \frac{r \, \chi \, \left[ \tilde{e}_+ \, {\rm erfc} \left( \alpha_+ \right)  + \tilde{e}_- \, {\rm erfc}\left( \alpha_- \right) \right]}{\chi \tilde{e}_+}  \right\rbrace \nonumber \\ 
 &= {\rm erfcinv} \left\lbrace r  \left[ {\rm erfc}\left( \alpha_+ \right)  + \frac{\tilde{e}_-}{\tilde{e}_+} \, {\rm erfc}\left( \alpha_- \right) \right] \right\rbrace  \nonumber \\
&= {\rm erfcinv} \left\lbrace  r  \left[ {\rm erfc}\left( \alpha_+ \right)  + \exp\left(-\frac{b\,c}{a} \right)  \, {\rm erfc}\left( \alpha_- \right) \right] \right\rbrace ,  \label{eq:leftinvcdf}
\end{eqnarray}
then, $y$ is given by $y = -(z - \alpha_+)/\sqrt{a}$. If it turns out that this $y$  fulfills $y > 0$, the other half of the cdf has to be inverted. Let
\begin{eqnarray}
 z &:= {\rm erfcinv} \left\lbrace  \left[ -\frac{r \, \mathcal{N}}{\chi} + \tilde{e}_+ {\rm erfc}\left( \alpha_+ \right)  + \tilde{e}_- {\rm erfc}\left( \alpha_- \right)\right] \tilde{e}^{-1}_- \right\rbrace \nonumber \\
&= {\rm erfcinv} \left\lbrace  (1-r) \left[ \exp\left(\frac{b\,c}{a} \right) {\rm erfc}\left( \alpha_+ \right) +  {\rm erfc}\left( \alpha_- \right)  \right] \right\rbrace . \label{eq:rightinvcdf}
\end{eqnarray}
Then, $y$ is given by $y = (z - \alpha_-)/\sqrt{a}$. \\
The complementary error function and its inverse are difficult to handle numerically, because there are no identities that allow to rescale or shift their evaluation to other intervals. Therefore, a robust numerical implementation of formulas \eref{eq:cdf},  \eref{eq:leftinvcdf} and  \eref{eq:rightinvcdf}  is rather involved. For the sake of a concise presentation, we present all details in \ref{sec:App_Imp}.

\subsection{MCMC Convergence Diagnostics} \label{subsubsec:MCMCConDia}

Assessing the efficiency of a sampling algorithm for a general purpose rather than a specific aim is a difficult task \cite{Li08}. Two types of \emph{convergence diagnostics} are usually applied: \emph{Qualitative diagnostics} rely on the visual inspection of some property of the chain $u_i,\, i=1,\ldots,K$. In contrast, \emph{quantitative diagnostics} try to compute characteristics that can be used to guide the sampling algorithm in an automated fashion. This should allow unexperienced users to perform ``black box'' Bayesian inference. Despite a lot of research on theses topics \cite{CoCa96,BrRo98,RoSa97,Th10}, no universal method is known. For our purpose, a qualitative autocorrelation analysis is appropriate. For a test function $g: \mathbb{R}^n \rightarrow \mathbb{R}^1$, the \emph{autocorrelation function} (\emph{acf}) $R:\{0,\ldots,K-1\} \rightarrow [-1,1]$ of the series $g_i := g(u_i),\, i=1,\ldots,K$ is given by:
\begin{eqnarray}
R(\tau) := \frac{1}{(K - \tau) \hat{\varrho}} \sum_{i=1}^{K -\tau} ( g_i - \hat{\mu} )   ( g_{i+\tau} - \hat{\mu} )  \\ 
\hat{\varrho} := \frac{1}{K} \sum_{i=1}^{K} ( g_i - \hat{\mu})^2, \quad \hat{\mu} := \frac{1}{K} \sum_{i=1}^{K} g_i
\end{eqnarray}
(Note that there are other possibilities to define $R$, but we need $R(0) = 1$). The value of $R(\tau)$ is referred to as the \emph{lag-$\tau$ autocorrelation} w.r.t. $g$. A fast decrease of the acf indicates that consecutive samples get mutually independent quite soon (if the $u_i$ would be independent, then, $R(\tau) = \delta_{(\tau,0)}$). For practical considerations, the decrease of autocorrelation w.r.t. to the raw number of samples drawn is not decisive if different samplers are compared. A method that has a slower decrease than others might still outperform them if it produces new samples considerably faster. In such situations, one would subsample the chain to get rid of highly correlated samples and to safe memory. Note that the notion of "one" sample is quite arbitrary anyway. In the SysGibbs sampler, one speaks of a "new" sample, if all components of $u$ are updated, in the MH-Si sampler one speaks of a "new" sample, if one component is updated. To address this, we will normally scale the acf by the computation time per sample $t_s$: $R^*(t) := R(t/t_s)$ for all $t=i \cdot t_s, i\in \{0,\ldots,K-1\}$, if we compare conceptually different sampling methods. $R^*(t) $ measures how fast a sampler can produce a certain loss in autocorrelation, which is of main interest for practical applications. However, while $R^*(t) $ is more decisive to compare different samplers, it relies on their concrete implementation\footnote{We implemented all samplers in Matlab and optimized them to yield the best possible performance. As mentioned in Section \ref{subsec:ConStr}, we originally indented to use the MH samplers in the scenario examined in Section \ref{subsec:EdgPreBayInv}. As their results were unsatisfactory even after a careful optimization of their implementation, we decided to develop the Gibbs samplers presented in this paper.}.\\
Normally, the test function $g$ is chosen with respect to the specific aim of inference. For instance, one could use the distance to the empirical mean of the whole chain if CM estimation is performed, or the projection onto a specific coordinate if that coordinate should be marginalized. Then, the rate of autocorrelation decrease is a measure of the efficiency of the chain for the specific inference aim. For our general purpose, we will test the ``worst case''. We project onto the direction of the largest variance, i.e., the first eigenvector $\nu_1$ of the covariance matrix $C$ of the posterior:
\begin{equation}
 g(u_i) := \langle \nu_1,u_i \rangle
\end{equation}
In general, the chain should have most problems to reduce the correlation of subsequent samples in this direction $\nu_1$. For each scenario we examine, the covariance matrix $C$ of the posterior is estimated from a long (sub-sampled) chain of the RnGibbs sampler, as this sampler will turn out to be the most reliable at a high performance. Note that this choice does not give an advantage to the RnGibbs sampler in the autocorrelation analysis but rather a disadvantage if the other samplers would have other directions of highest variance. We checked that this is not the case in a test scenario we examined in preliminary studies. \\
Other possible MCMC convergence diagnostic plots that are commonly used are plots of $\log[p(u_i|m)]$ or of single components $(u_i)_j$. Such plots are good to detect possible \emph{multimodality} of the posterior and to determine a sufficient number of burn in steps $K_0$. Multimodality is not an issue in our case, as the posterior is log-concave (the energy $-\log[p(u|m)]$ is convex). The burn-in length is an important factor for the practicability of the algorithms (and we will address this issue in our studies) but it is a difficult measure for a fair and definite comparison of the sampling methods. First, it crucially relies on the initialization of the chain, so one would have to compare all methods for various common initialization strategies, which is not really feasible and too application specific. Second, for the Metropolis-Hastings schemes, an adaptation of the sampling parameters to is usually carried out in the burn-in phase with the aim to optimize the performance of the chain in the real run. We will introduce this topic in Section \ref{para:ChoicePara}. The consequence is that $K_0$ also depends on the adaptation scheme, which renders the problem of a meaningful comparison even worse.


\section{Results}\label{sec:Results}

In this section, we compare the sampling algorithms for two scenarios: Edge-preserving, TV-based image deblurring in 1D and impulse prior based image deblurring in 2D. All algorithms have been implemented in Matlab and have been optimized to the best possible performance. All results have been computed on the same CPU architecture limiting Matlab to a single computational thread, i.e., to use a single CPU core with 2.80GHz (parallelization is discussed in Section \ref{sec:Discussion}). We paid special attention that the computation times are as comparable as possible.

\subsection{Edge-Preserving Bayesian Inversion in 1D} \label{subsec:EdgPreBayInv}
A popular case of L1-type priors arises from \emph{edge-preserving} image reconstruction. The task is to reconstruct a spatially distributed intensity image that is known to consist of piecewise homogeneous parts with sharp edges from indirect, noisy measurements (e.g., the recovery of the body's organs and their boundaries from X-ray computed tomography data \cite{KoSiJaKaKoLaPiSo03,SiKoJaKaKoLaPiSo03,KaSo05}). Using Gaussian, i.e., L2 -type priors smooths the image edges in such situations. In contrast, \emph{total variation} (\emph{TV}) priors, which rely on the L1 norm of the first spatial derivatives, are able to retain them \cite{RuOsFa92,KaSo05,Lo08,BuOs12}. The use of TV priors in Bayesian inference has led to interesting theoretical questions. It was discovered that it is not possible to formulate the conventional TV prior in a \emph{discretization invariant} way \cite{LaSi04,LaSaSi09}, i.e., that the posterior converges to a well defined limit probability density when the level of discretization is increased while reflecting the a priori information of edge-preservation at all levels of discretization. If the TV prior is formulated such that it converges, it converges to a Gaussian smoothness prior, and, thus, the edge-preservation property is lost. This motivated research on whether and how it is possible to formulate edge-preservation as a priori information in a consistent, discretization invariant way in the Bayesian framework. Recently, \emph{Besov space priors} have been proposed, which rely on a weighted L1 norm of wavelet basis coefficients \cite{LaSaSi09,KoLaNiSi12}. Such priors are L1-type priors with an invertible $D$, thus, posterior sampling by means of our Gibbs sampling algorithms can be performed with ease. In addition, modifications of the standard TV prior \cite{Co11} and \emph{hierarchical Bayesian models} \cite{CaSo07,CaSo08,HeLa09,He10,BaCaSo10} have been proposed for discretization invariant edge-preserving image reconstruction as well. \\
To address the problem of discretization invariance of a prior, one can, e.g., study the convergence of the corresponding CM estimate for $n \rightarrow \infty$. The problems of using MH-based samplers for CM estimation in high dimensions have already been noticed in \cite{LaSi04,KoLaNiSi12} and our research on alternative samplers has been motivated by these problems as well. 

\subsubsection{Setting} \label{subsubsec:Set}

\begin{figure}[hbt]
   \centering
\subfigure[][ \label{subfig:Real}]{\includegraphics[height=2.65cm]{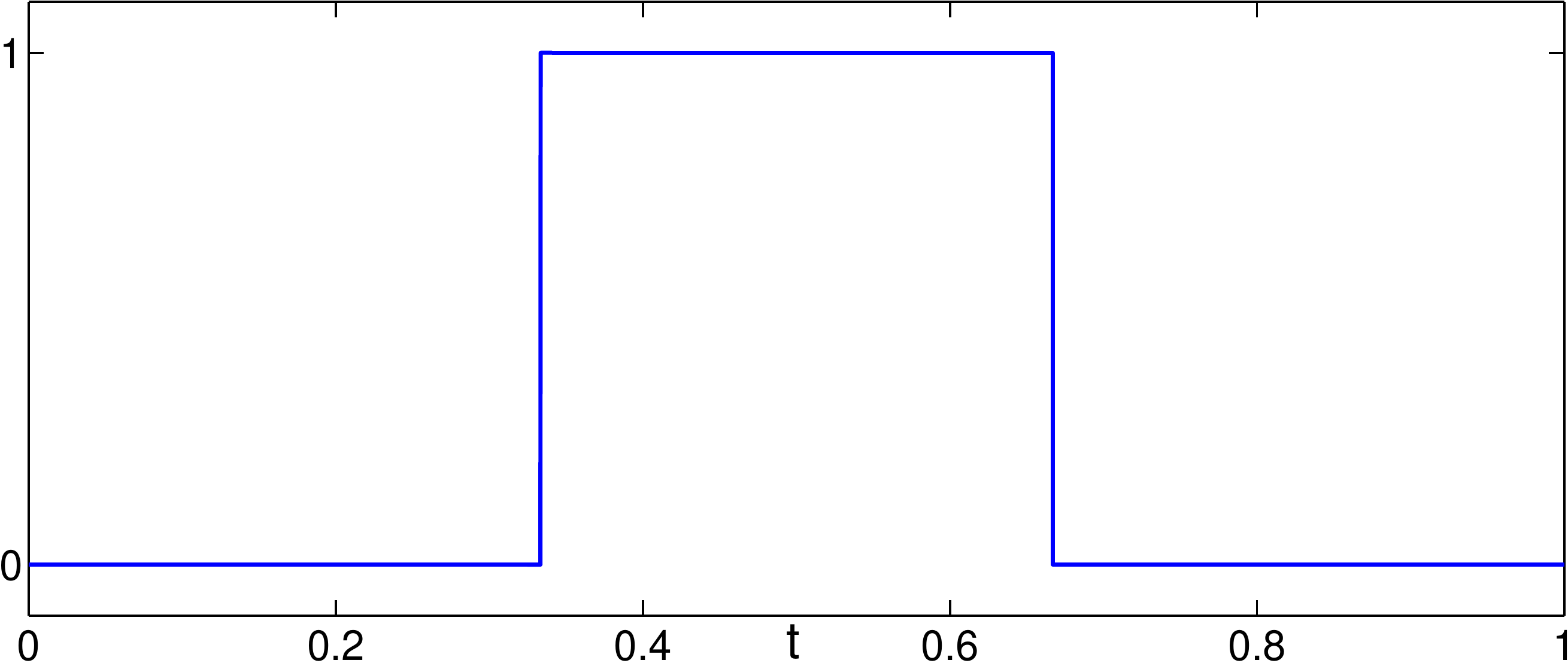}}
\subfigure[][  \label{subfig:Data}]{\includegraphics[height=2.65cm]{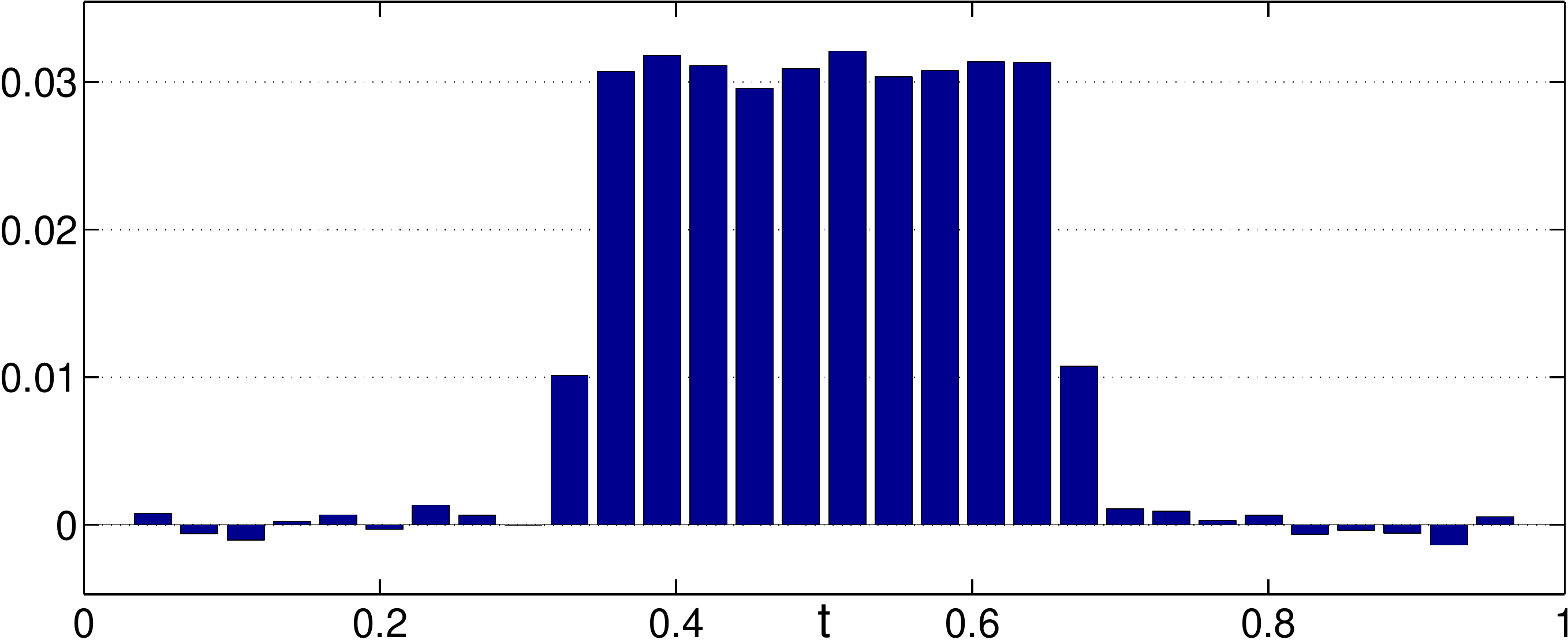}}
\caption{Left: The unknown function $\tilde{u}(t)$. Right: The measurement data $m$. }
   \label{fig:RealData}
\end{figure}

We rely on the setting used in \cite{LaSi04}. The motivation is to mimic a measurement made by a \emph{charge coupled device} (\emph{CCD}) used in digital cameras or medical imaging devices. These devices integrate the amount of light illuminating a certain pixel over a certain period of time. In the continuous model setting, we represent the unknown light intensity by a positive function $\tilde{u}:[0,1] \rightarrow \mathbb{R}$, and the $k$ pixels of the CCD device as a equidistant division of the subinterval $[{\small \frac{1}{k+2}},{\small \frac{k+1}{k+2}}] \subset [0,1]$, i.e., the $j$-th pixel is represented by the interval $[{\small \frac{j}{k+2}},{\small \frac{j+1}{k+2}}] $. The measurement at the $j^{th}$ pixel is then given by:
\begin{equation}
 m_j  = \int_{\frac{j}{k+2}}^{\frac{j+1}{k+2}} \tilde{u}(t) \rmd t + \varepsilon_j \label{eq:contfwd}
\end{equation}
For convenience, we will choose $k = 2^{L_m} - 2$ and $L_m = 5$. For discretizing $\tilde{u}$, we choose the grid $\{{\small \frac{1}{n+1}},\ldots,{\small \frac{n}{n+1}}\} \subset [0,1]$,  and let $n = 2^{L_u} - 1$ with $L_u > L_m$. The discretization of the forward mapping implied by \eref{eq:contfwd} in terms of the $k \times n$ matrix $A$ can then be implemented by the trapezoidal quadrature rule. The $j^{th}$ row of $A$ is given by
\begin{equation}
A_{(j,\cdot)} := [\underbrace{0,0,\ldots,0,}_{j2^{(L_u-L_m)} - 1} \frac{1}{2}h ,\underbrace{h,h,\ldots,h}_{2^{(L_u-L_m)} - 1},\frac{1}{2}h,0,0,\ldots,0],
\end{equation} 
where $h := {\small \frac{1}{n+1}}$ defines the grid size. The discrete TV prior with Neumann boundary conditions in our situation is given by:
\begin{eqnarray}
p(u) \propto \exp \left( - \lambda_n \sum_{i=1}^{n-1} |u_{i+1} - u_i| \right) = \exp \left( - \lambda_n | D \, u | \right) \\
\end{eqnarray}
where $D \in \mathbb{R}^{(n-1) \times n}$ is given by $(D u)_i :=  u_{i+1} - u_i$, $i=1,\ldots,n-1$. We indexed $\lambda_n$ by $n$ to stress that we may choose it depending on the discretization level.\\
For the Gibbs sampler (cf. Sections \ref{subsec:GenSetBayFor} and \ref{subsec:ImpGibbsL1}), we note that $l = {\rm rank}(D) = n-1$, and $v_1,\ldots,v_{n-1} \in \mathbb{R}^n$ are given by step functions: $(v_i)_j =  \mathbf{1}_{\{ j > i \}}$. These are completed to a basis of $\mathbb{R}^n$ by $(v_n)_j = 1 \; \forall j$. If we reorder them and define $V := [v_n,v_1,\ldots,v_{n-1}]$, we can write $V$ as
\begin{equation}
 V_{(i,j)} = \cases{1 & if $i \geqslant j$\\
	    0 & else}
\end{equation}
The unknown function $\tilde{u}$ we actually use is the indicator function on $[{\small \frac{1}{3}},{\small \frac{2}{3}}]$, see Figure \ref{subfig:Real}. Measurement data $m$ is generated using formula \eref{eq:contfwd}, see Figure \ref{subfig:Data}. The standard deviation of the measurement noise $\sigma$ is 0.001.
We will examine different combinations of $n$ and $\lambda_n$:
\begin{enumerate}
 \item[A:] $n = 63$ in combination with $\lambda_n = 100$, $200$ and $400$, respectively. Here, we focus on increasing the impact of the prior. The posterior will become less Gaussian because the weight of the L1-type TV prior is increased.
 \item[B:] $n = 2^{L_u} - 1$ for $L_u = 7,8,\ldots$ with $\lambda_n = 25 \cdot \sqrt{n+1}$. With this scaling of $\lambda_n$, the posterior $p(u|m)$ converges for $n \longrightarrow \infty$, but the edge-preserving property of the TV prior is lost, see \cite{LaSi04} for details. The CM estimate $u_{\rm CM}$ converges to a smooth limit function, which will facilitate the visual validation of the results of the different MCMC methods.
\end{enumerate}

\subsubsection{Preliminaries} \label{subsubsec:Pre}

\paragraph{Choice of Parameters:} \label{para:ChoicePara}
For the MH schemes, the proper tuning of $\kappa$ is essential. If it is very small, the proposals will always be accepted since the distribution to sample from is continuous. However, in return the exploration of the sampling space is slow. On the contrary, if is too large, the differences in probability will be huge because the distribution is log-concave and new proposals will hardly be accepted. A good overview on this topic is given in \cite{RoRo01,NeRo06}. The remarkable result is that in high dimensions, having a total acceptance rate of new proposals of about 0.234 leads to an optimal efficiency independent of the distribution to sample from. Furthermore, this optimal efficiency hardly drops in the range between 0.1 and 0.4 of acceptance rate. This yields an easy to implement rule to tune $\kappa$: One could find the optimal $\kappa$ in a preliminary MH-MCMC run and initialize the real MH-MCMC run with it. However, it turns out that this $\kappa$ is only optimal once the chain has reached the main support of the distribution while it can hinder the chain from ever getting there (the burn-in length increases dramatically, see Section \ref{subsubsec:MCMCConDia}). For these reasons, on-line adaptation of $\kappa$ is usually used. The empirical acceptance rate is monitored, and $\kappa$ is increased if it is too high while $\kappa$ is decreased if it is too low. The scheme we use is that every $10\,000$ samples, the empirical acceptance rate is computed and if it is above 0.35, $\kappa$ is multiplied by 1.2 while it is multiplied by 0.8 if it is below 0.15. In theory, the resulting chain will then not be a Markov chain anymore (but it is still ergodic). However, in practice, using this scheme, $\kappa$ hardly ever changes once the burn-in time is over and so the real chain is not affected.\\
For MH-Ncom, we have to choose $n^*$, i.e., the number of components that are updated in one step. We choose $n^* = \lfloor n^{7/12} \rfloor$, which roughly corresponds to the values used in \cite{LaSi04}.

\begin{figure}[hbt]
  \centering
\includegraphics[width = \textwidth]{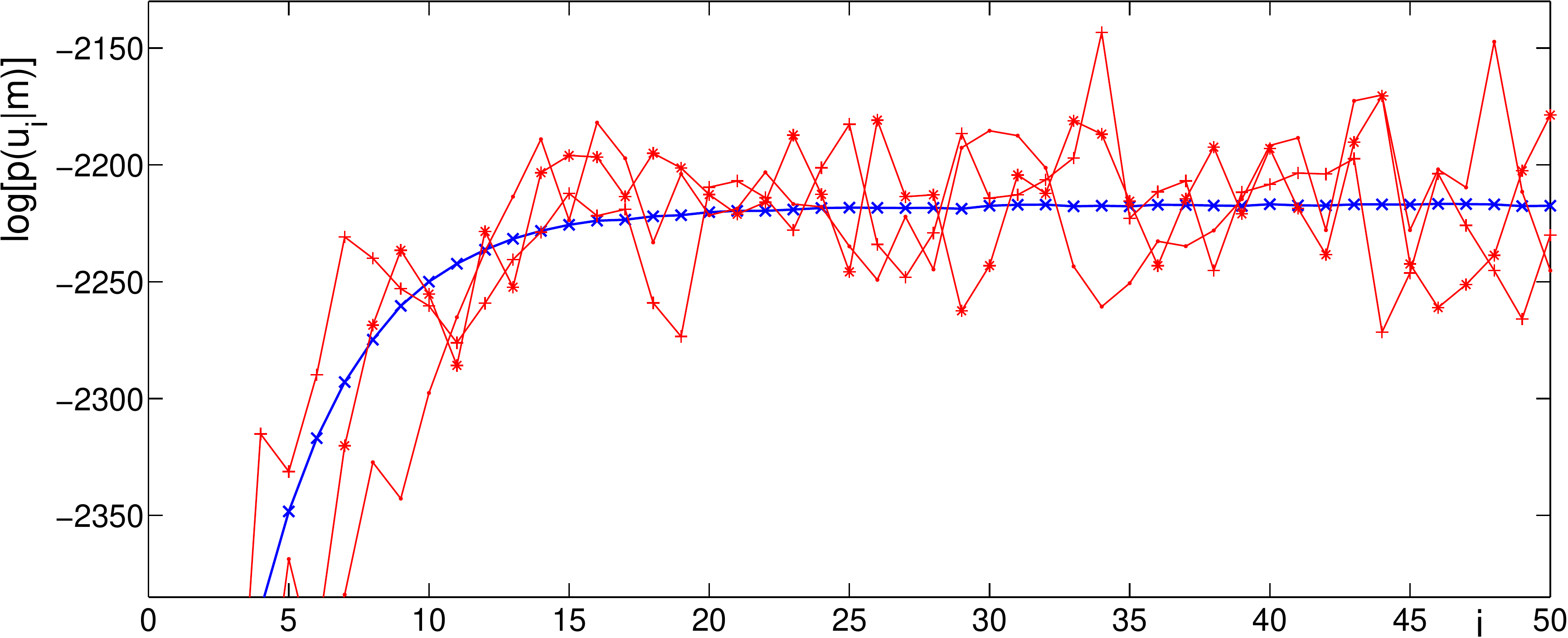}
\caption{Plots of  $\log[p(u_i|m)]$ for $n$ = 1023, $\lambda$ = 800, using the RnGibbs sampler.  Red \mbox{$\bullet$}, \mbox{$\ast$}, \textbf{+}: three independent realizations. Blue \mbox{$\times$}: the average of 5000 independent realizations.
   \label{fig:burnin}}
\end{figure}

\paragraph{Burn-in Times:} \label{para:BurnIn}
As noted in Section \ref{subsubsec:MCMCConDia}, the sufficient amount of burn-in steps that have to be drawn can be deduced from observing $\log[p(u_i|m)]$. Once it oscillates around a constant value, the stationary part of the distribution is reached. Averaging  $\log[p(u_i|m)]$ over a large number of independent chains that all started at the same initialization ($u \equiv 0$ in our case) removes the oscillations and allows to determine $K_0$ in an easy fashion. See Figure \ref{fig:burnin} for an example of such a plot. In Table \ref{tbl:burn-in}, the burn-in steps $K_0$ and the corresponding computation times $t_0$ are listed for the combinations of  $n$ and $\lambda$ that are examined in detail. It gives a first impression of how the methods scale with $n$ and $\lambda_n$, but as noted in Section \ref{subsubsec:MCMCConDia} it does not allow for a fair and detailed comparison.

\begin{table}
\caption{\label{tbl:burn-in} Necessary burn-in steps $K_0$ and computation time $t_0$ in seconds for each method and combination of $n$ and $\lambda$ when starting from $u \equiv 0$. The values were found by observing the log posterior probability averaged over a large number of independent chains. Listed as $(K_0 , t_0)$ }
\begin{indented}
\lineup
\item[] \begin{tabular}{@{}*{4}{l}}
\br
 &\centre{3}{Model parameters ($n$,$\lambda$)}\\
\ns
Method &&&\\
& (63,100) & (63,200) & (63,400)  \\
\mr
MH-Iso 	          & (4e5,1.8e1)    & (4e5,1.9e1)   & (5e5,2.3e1)   \\
MH-Ncom          & (4e5,2.3e1)    & (4e5,2.5e1)   & (5e5,2.9e1)   \\
MH-Si 		  & (5e5,2.8e1)    & (5e5,3.0e1)   & (6e5,3.4e1)   \\
RnGibbs 	          & (200,0.5e0)    & (200,0.5e0)   & (200,0.4e0)   \\
RnGibbsO3 	  & (200,0.9e0)    & (200,1.0e0)   & (200,0.9e0)   \\
RnGibbsO7 	  & (200,1.0e0)    & (200,1.0e0)   & (200,0.9e0)   \\
SysGibbs 	          & (400,1.0e0)    & (500,1.3e0)   & (500,1.3e0)   \\
SysGibbsO3  	  & (200,0.9e0)    & (500,2.2e0)   & (500,2.2e0)   \\
SysGibbsO7  	  & (200,0.9e0)    & (500,2.3e0)   & (500,2.2e0)  \\
\br
\end{tabular}
\item[] \begin{tabular}{@{}*{5}{l}}
\br
 &\centre{4}{Model parameters ($n$,$\lambda$)}\\
\ns
Method &&&\\
&  (127,280) & (255,400) & (511,560) & (1023,800) \\
\mr
MH-Iso 	          & (7e5,3.4e1) & (4e6,2.1e2)        & (3e7,1.9e3)  & (2e8,1.7e4) \\
MH-Ncom          & (7e5,4.2e1) & (4e6,2.6e2)           & (3e7,2.5e3)  & (2e8,2.3e4 )\\
MH-Si 		  & (8e5,4.4e1) & (4e6,2.3e2)      & (3e7,1.9e3)  & (2e8,1.5e4) \\
RnGibbs 	          & (80,0.4e0)   &   (50,0.4e0)      & (30,0.5e0)    & (20,0.6e0) \\
RnGibbsO3 	  & (80,0.7e0)   &   (50,0.9e0)      & (30,1.1e0)    & (20,1.4e0) \\
RnGibbsO7 	  & (80,0.7e0)   &   (50,0.9e0)      & (30,1.1e0)    & (20,1.5e0) \\
SysGibbs 	          & (150,0.7e0) & (100,0.9e0)  & (150,2.9e0)  & (150,5.3e0)  \\
SysGibbsO3  	  & (150,1.4e0) & (150,2.6e0)      & (150,5.2e0)  & (150,1.0e1)  \\
SysGibbsO7       & (150,1.4e0) & (150,2.7e0)      & (200,6.8e0)  & (200,1.4e1)  \\
\br
\end{tabular}
\end{indented}
\end{table}

\subsubsection{General Autocorrelation Analysis} \label{subsubsec:GenAutAna}
\begin{figure}[hbt]
   \centering
\subfigure[][MH-based samplers \label{fig:acf_n63_MH}]{\includegraphics[height=0.47 \textwidth]{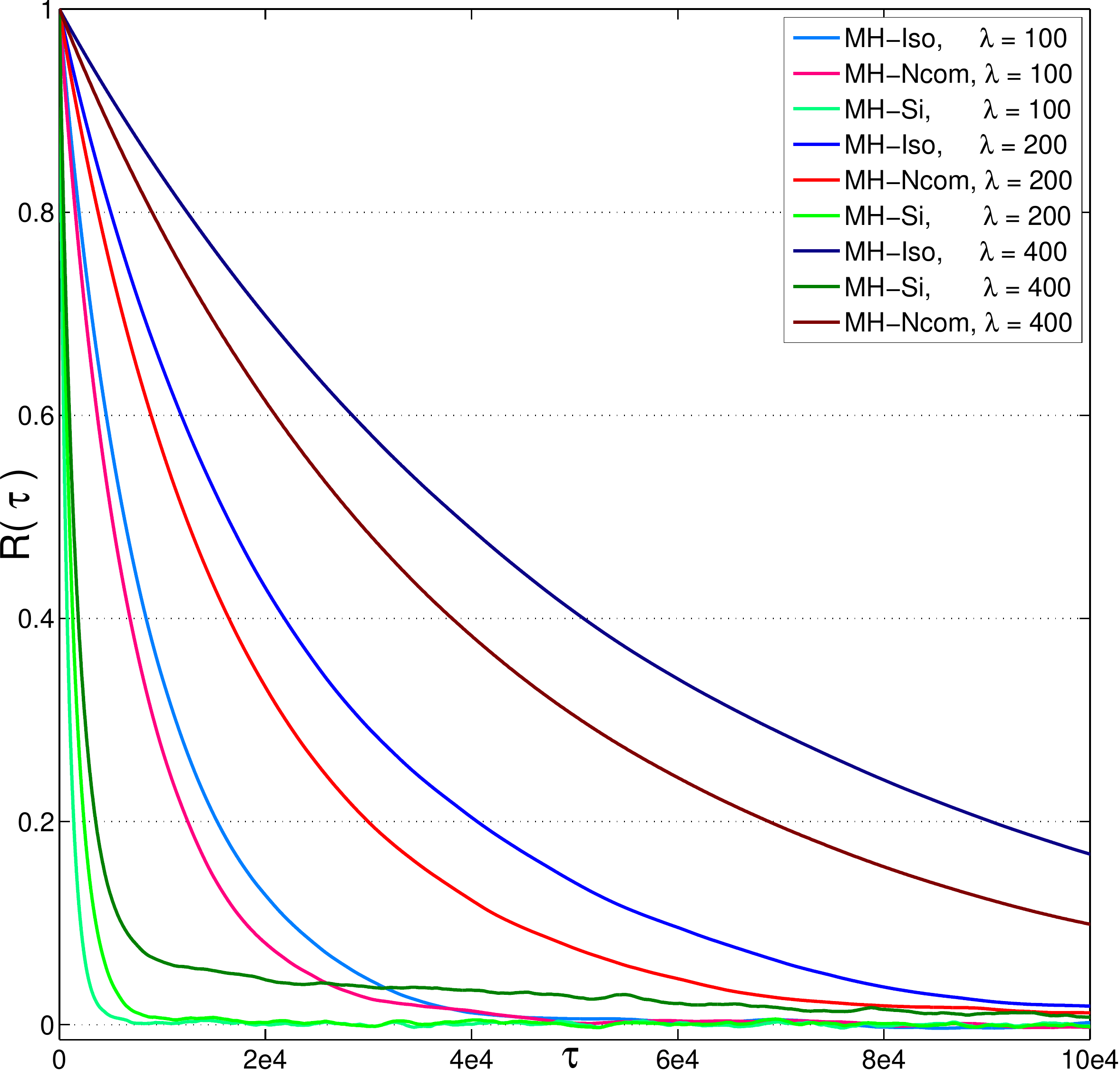}}
\subfigure[][Gibbs-based samplers \label{fig:acf_n63_Gibbs}]{\includegraphics[height=0.47 \textwidth]{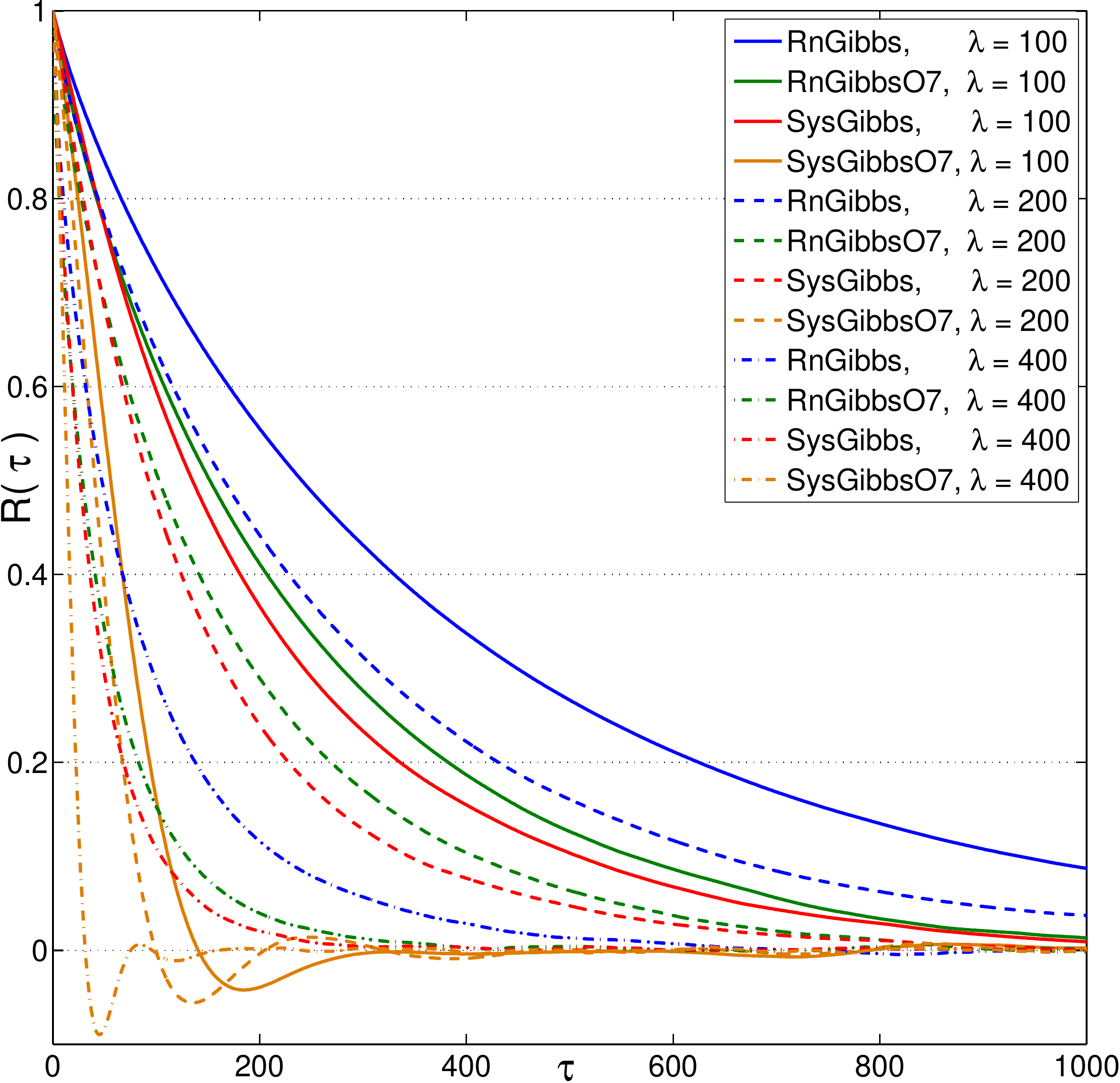}}
\caption{Autocorrelation plots $R(\tau)$ for $n = 63$, $\lambda_n = 100$, $200$ and $400$. \label{fig:acf_n63}}
\end{figure}

\begin{figure}[hbt]
   \centering
\subfigure[][MH vs. un-overrelaxed Gibbs sampler \label{fig:tacf_n63_MHvsGibbs}]{\includegraphics[width=0.49 \textwidth]{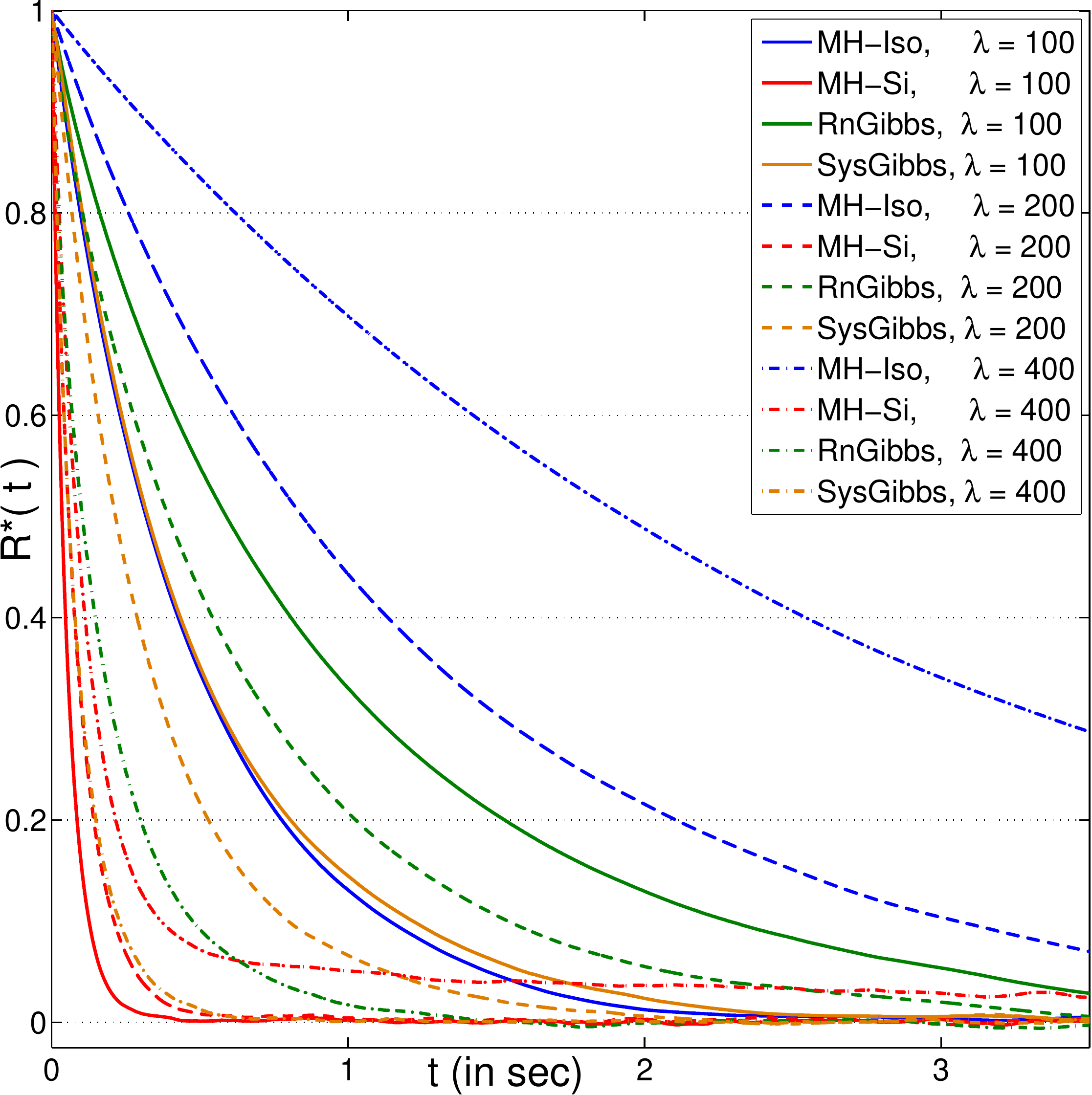}}
\subfigure[][Normal vs. overrelaxed Gibbs samplers \label{fig:tacf_n63_Gibbs}]{\includegraphics[width=0.49 \textwidth]{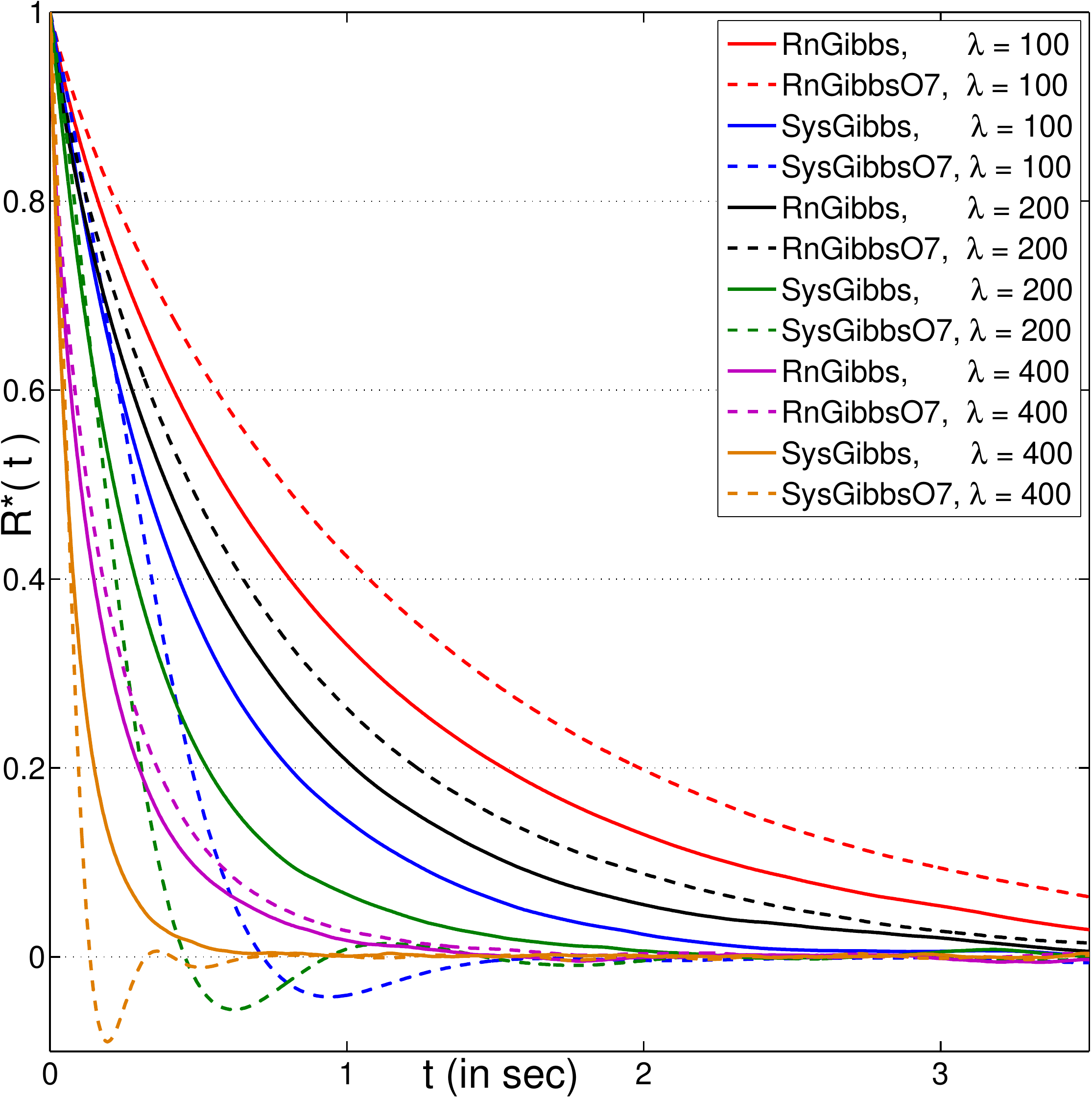}}
\caption{Temporal autocorrelation plots $R^*(t)$ for $n = 63$, $\lambda_n = 100$, $200$ and $400$. \label{fig:tacf_n63}}
\end{figure}

\begin{figure}[hbt]
   \centering
\includegraphics[width = \textwidth]{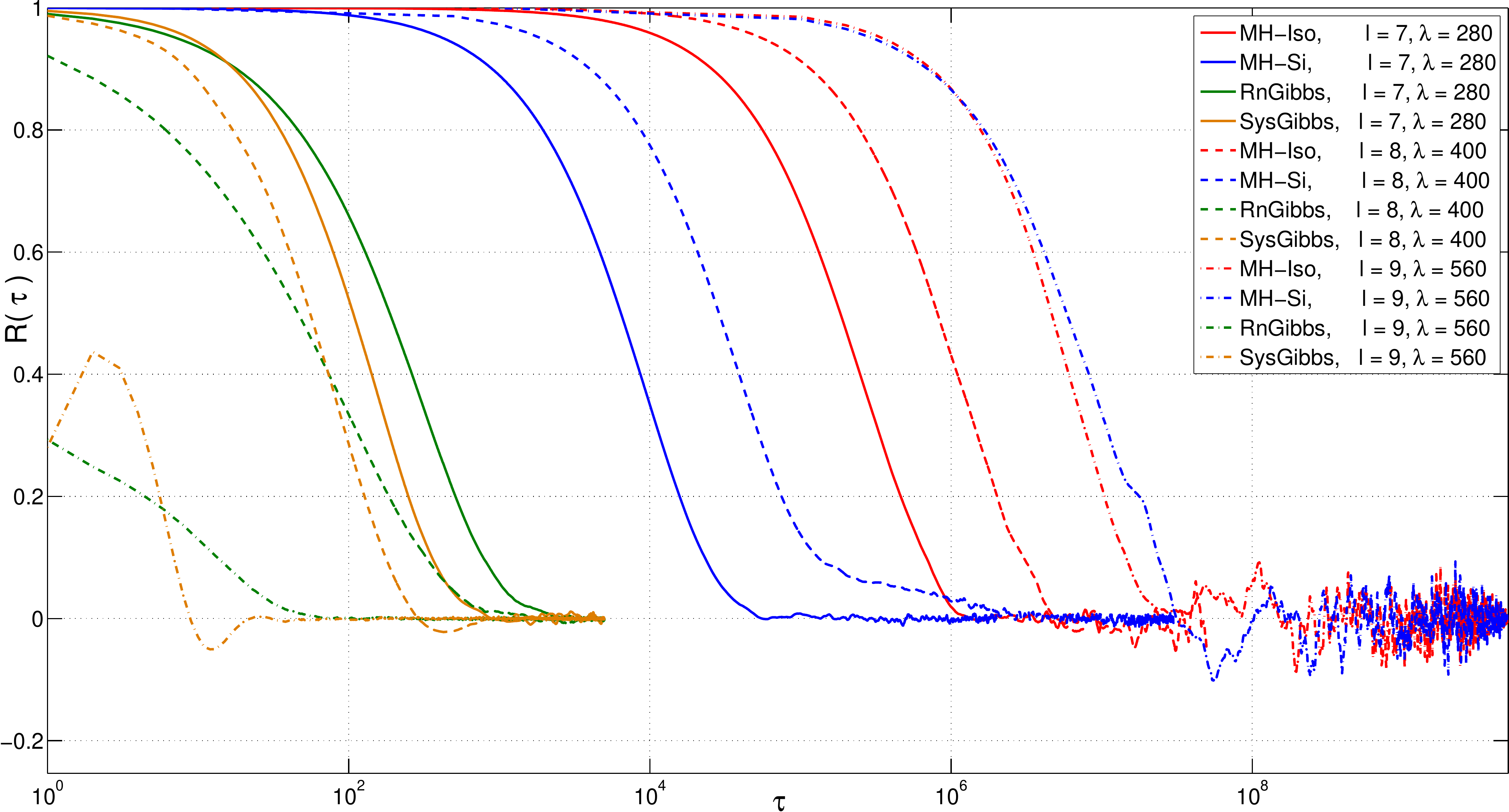}
\caption{Autocorrelation plots $R(\tau)$ for varying $n$ and $\lambda$ for MH and un-overrelaxed Gibbs sampler. Note that the $\tau$ axis starts at $\tau = 1$ and is scaled logarithmically. \label{fig:acf_nVar_MHvsGibbs}}
\end{figure}

\begin{figure}[hbt]
   \centering
\includegraphics[width = \textwidth]{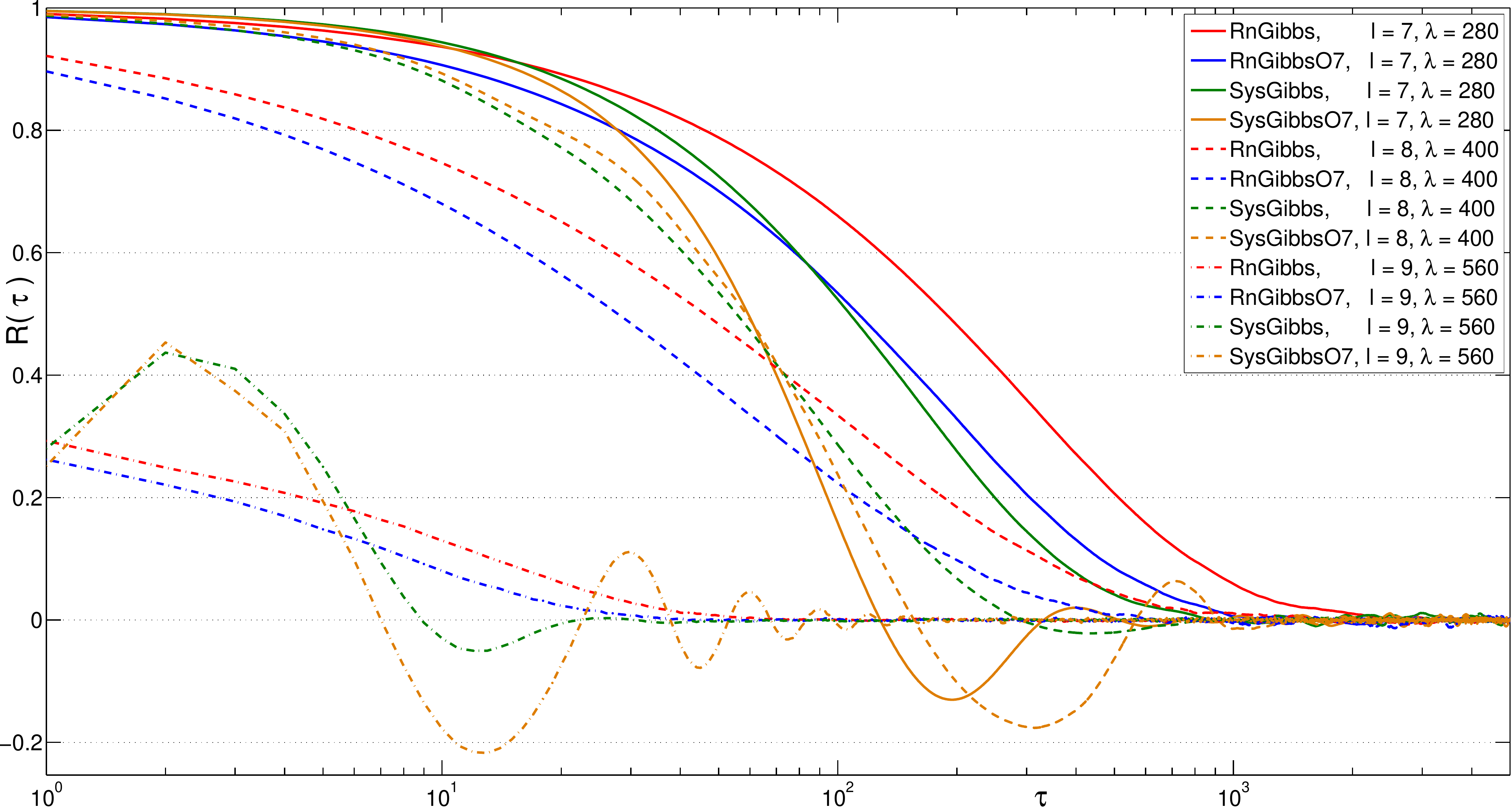}
\caption{Autocorrelation plots $R(\tau)$ for varying $n$ and $\lambda$ for overrelaxed Gibbs sampler. Note that the $\tau$ axis starts at $\tau = 1$ and is scaled logarithmically. \label{fig:acf_nVar_Gibbs}}
\end{figure}

\begin{figure}[hbt]
   \centering
\includegraphics[width = \textwidth]{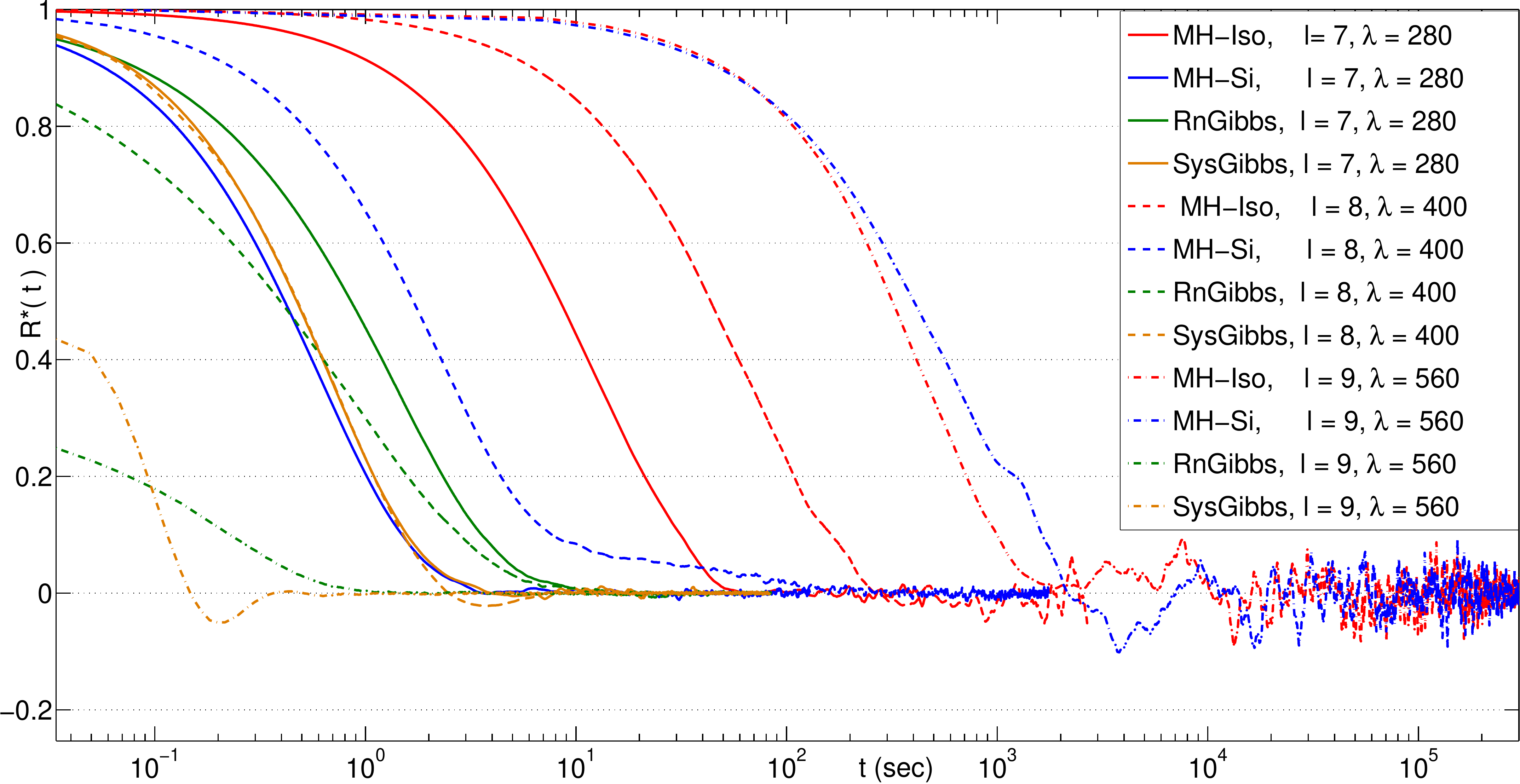}
\caption{Temporal autocorrelation plots $R^*(t)$ for varying $n$ and $\lambda$ for MH un-overrelaxed Gibbs sampler. Note that the $t$ axis starts at the smallest $t_s$ and is scaled logarithmically. \label{fig:tacf_nVar_MHvsGibbs}}
\end{figure}

As explained in Section \ref{subsubsec:MCMCConDia}, we will rely on autocorrelation plots for the projection of the samples onto the direction of maximal covariance as qualitative measures of the efficiency of the sampling algorithms. Figure \ref{fig:acf_n63} shows the autocorrelation plots $R(\tau)$ for $n = 63$ and varying $\lambda$. In Figure \ref{fig:tacf_n63}, the corresponding temporal autocorrelation plots $R^*(t)$ are shown. The comparison is split up into MH-based vs. normal Gibbs samplers and normal vs. overrelaxed Gibbs samplers to reduce the number of plots shown in one figure. The plots for RnGibbsO3, SysGibbsO3 and MH-Ncom were omitted for the same reason: The plots for RnGibbsO3 and SysGibbsO3 lie between the plots of RnGibbs and RnGibbsO7 and SysGibbs and SysGibbsO7, respectively. The plots of MH-Ncomp look similar to the ones of MH-Iso and lie between MH-Iso and MH-Si. Figures \ref{fig:acf_nVar_MHvsGibbs} and  \ref{fig:acf_nVar_Gibbs} show the autocorrelation plots $R(\tau)$ for varying $n$ and $\lambda_n = 25 \cdot \sqrt{n+1}$. In Figure \ref{fig:tacf_nVar_MHvsGibbs}, the temporal autocorrelation plots $R^*(t)$ corresponding to Figure \ref{fig:acf_nVar_MHvsGibbs} are shown. The plots for RnGibbsO3, SysGibbsO3 and MH-Ncom are, again, omitted. In addition, the plots for $l = 10$ are not shown as the trends of the autocorrelation functions for growing $n$ are already clearly visible.\\
Table \ref{tbl:lag1} lists the lag $\tau_{0.01}$ for which the autocorrelation $R(\tau)$ drops below 1 \% for the first time and the corresponding computation time $t_{0.01} = \tau_{0.01} \cdot t_s$.

\begin{table}
\caption{\label{tbl:lag1} Lag of 1\% auto correlation for each method and combination of $n$ and $\lambda$ in terms of (samples,computation time)}
\begin{indented}
\lineup
\item[] \begin{tabular}{lrrr}
\br
 &\centre{3}{Model parameters ($n$,$\lambda$)}\\
\ns
Method &&&\\
& (63,100) & (63,200) & (63,400)  \\
\mr
MH-Iso 	          & (4.1e4,2.1e0)   & (1.2e5,6.2e0)   & (2.1e5,1.1e1)\\
MH-Ncom          & (4.2e4,2.6e0)   & (1.0e5,6.4e0)   & (2.1e5 ,1.0e1) \\
MH-Si 		  & (4.7e3,0.3e0)     & (7998,0.5e0)       & (8.8e4,5.4e0) \\
RnGibbs     	  & (1685,4.1e0)     & (1402,3.3e0)       & (561,1.2e0) \\
RnGibbsO3 	  & (1239,5.8e0)     & (983,4.6e0)           & (395,1.8e0) \\
RnGibbsO7 	  & (1056,5.5e0)     & (811,3.7e0)         & (318,1.4e0) \\
SysGibbs   	  & (985,2.4e0)       & (810,1.9e0)           & (242,0.5e0) \\
SysGibbsO3       & (412,1.9e0)         & (242,1.1e0)         & (76,0.3e0) \\
SysGibbsO7       & (137,0.7e0)         & (97,0.4e0)         & (31,0.1e0) \\
\br
\end{tabular}

\item[] \begin{tabular}{lrrrr}

 &\centre{4}{Model parameters ($n$,$\lambda$)}\\
\ns
Method &&&\\
&  (127,280) & (255,400) & (511,560) & (1023,800) \\
\mr
MH-Iso 	             & (1.1e6,5.0e1) & (4.6e6,2.5e2) & (2.7e7,1.9e3) & (1.3e8,1.3e4) \\
MH-Ncom              & (9.4e5,5.4e1)    & (3.2e6,2.2e2) & (2.1e7,1.8e3) & (1.8e8,2.3e4) \\
MH-Si 		     & (4.8e4,3.1e0)   &  (2.1e6,1.2e2)& (3.1e7,2.1e3) & (2.9e8,2.5e4) \\
RnGibbs 	     & (2017,9.2e0)    & (1014,8.7e0)    & (46,0.8e0)     & (39,1.3e0) \\
RnGibbsO3 	     & (1006,8.9e0)        & (1052,2.0e1)    & (31,1.1e0)   & (29,2.1e0) \\
RnGibbsO7 	     & (953,8.7e0)       & (473,8.7e0)       & (28,0.9e0) & (24,1.8e0) \\
SysGibbs 	     & (770,3.4e0)       & (270,2.3e0)       & (9,0.1e0)     & (12,0.4e0)\\
SysGibbsO3           & (230,2.0e0)          & (165,2.9e0)       & (8,0.3e0)   &  (7,0.5e0) \\
SysGibbsO7           & (126,1.2e0)       & (153,2.8e0)       & (7,0.2e0)   & (6,0.4e0) \\
\br
\end{tabular}

\end{indented}
\end{table}

\subsubsection{Visual Results} \label{subsubsec:VisRes1D}
\begin{figure}[hbt]
  \centering
\subfigure[][MH-Iso\label{subfig:visu_MH-Iso}]{\includegraphics[width=0.3\textwidth]{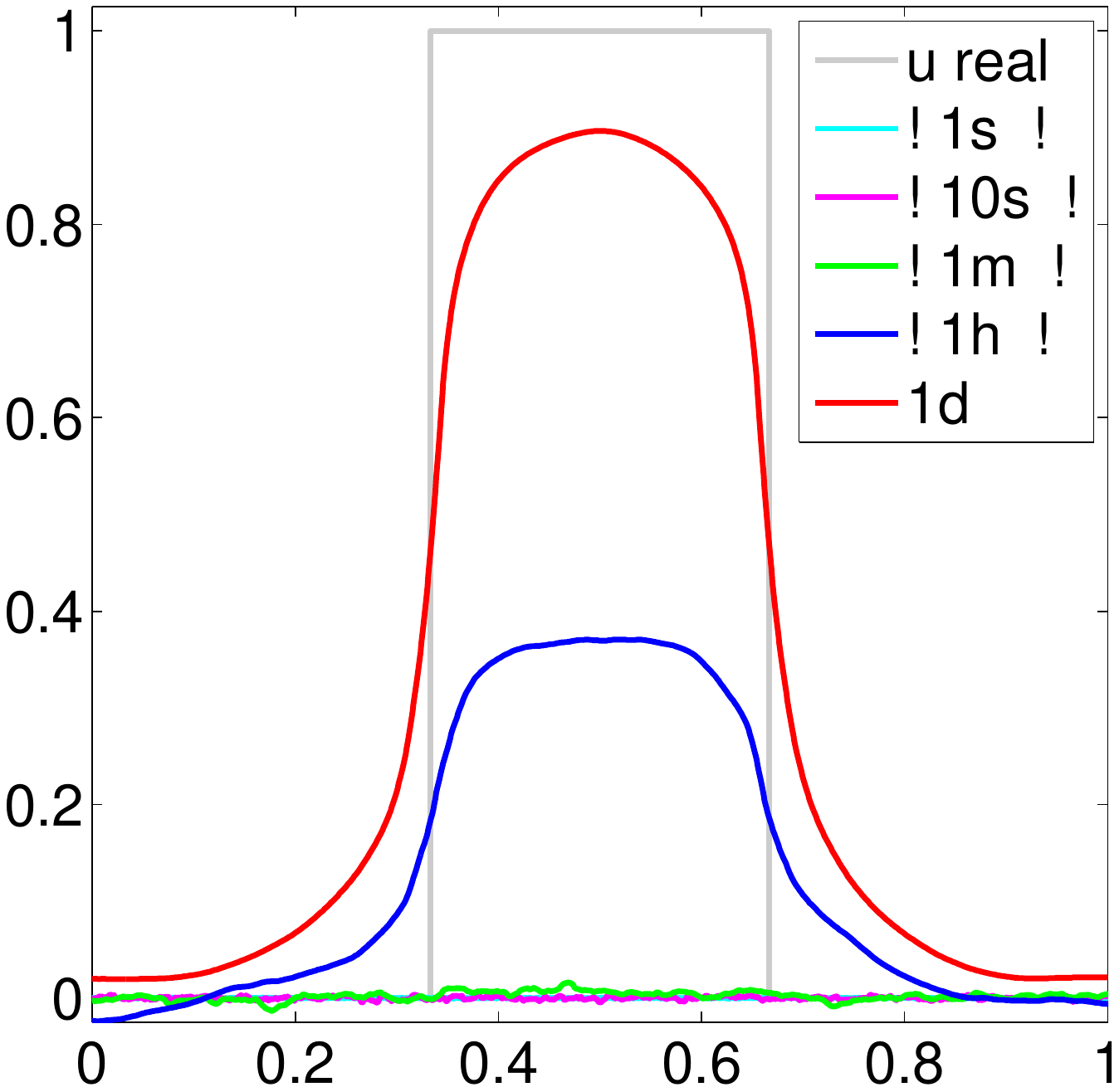}}
\subfigure[][MH-Ncom\label{subfig:visu_MH-Ncom}]{\includegraphics[width=0.3\textwidth]{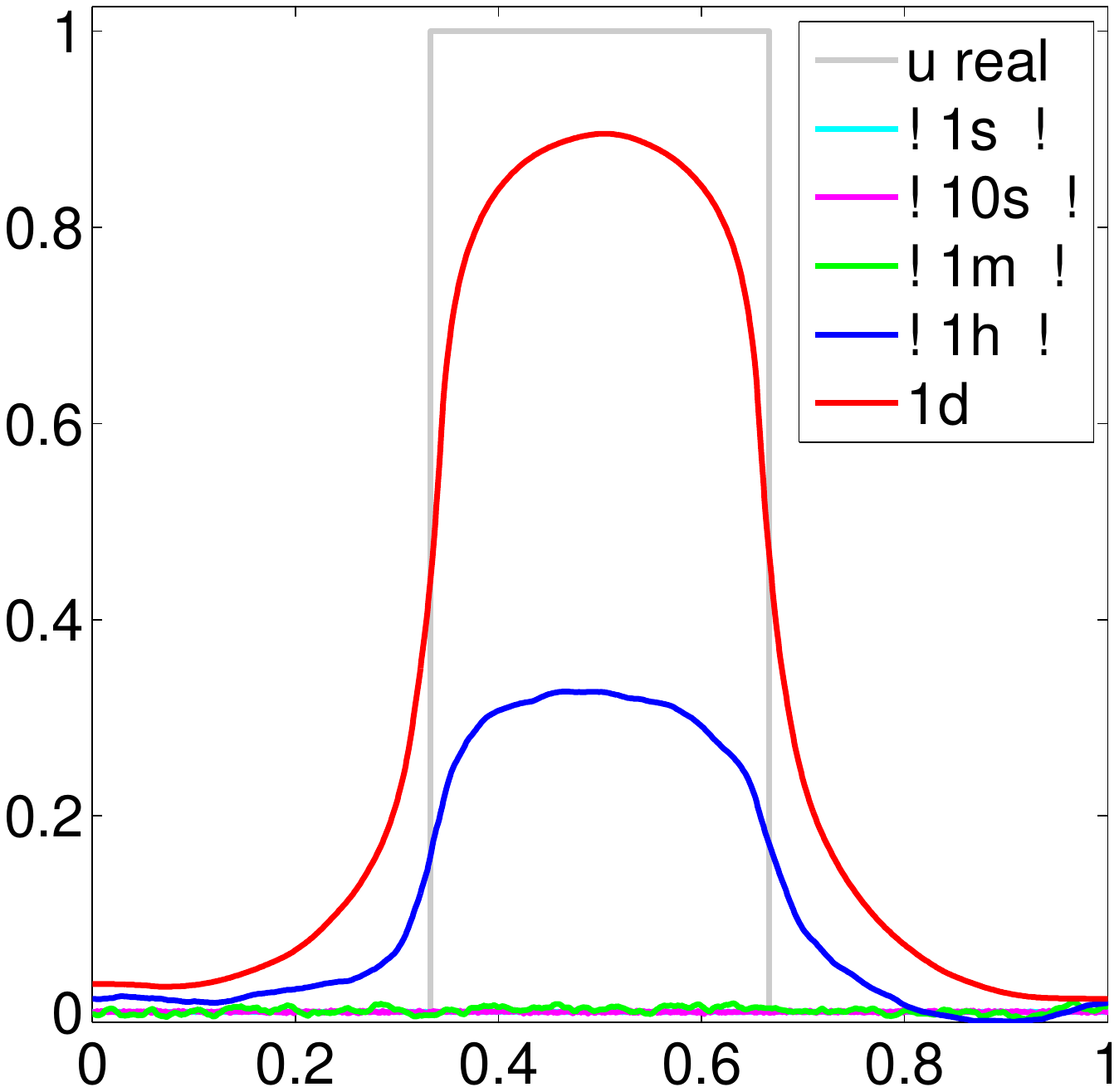}}
\subfigure[][MH-Si\label{subfig:visu_MH-Si}]{\includegraphics[width=0.3\textwidth]{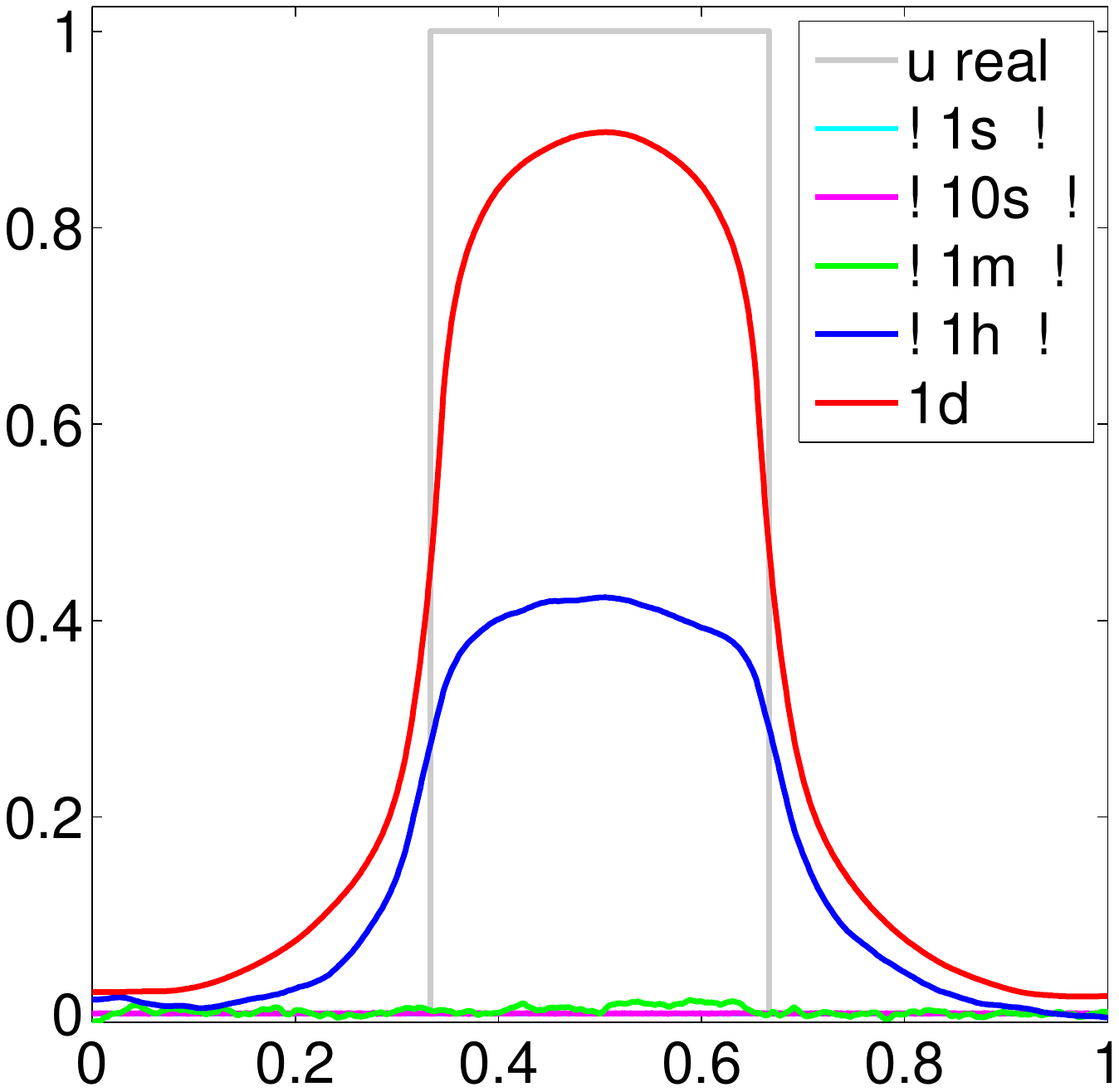}}\\
\subfigure[][RnGibbs\label{subfig:visu_RnGibbs}]{\includegraphics[width=0.3\textwidth]{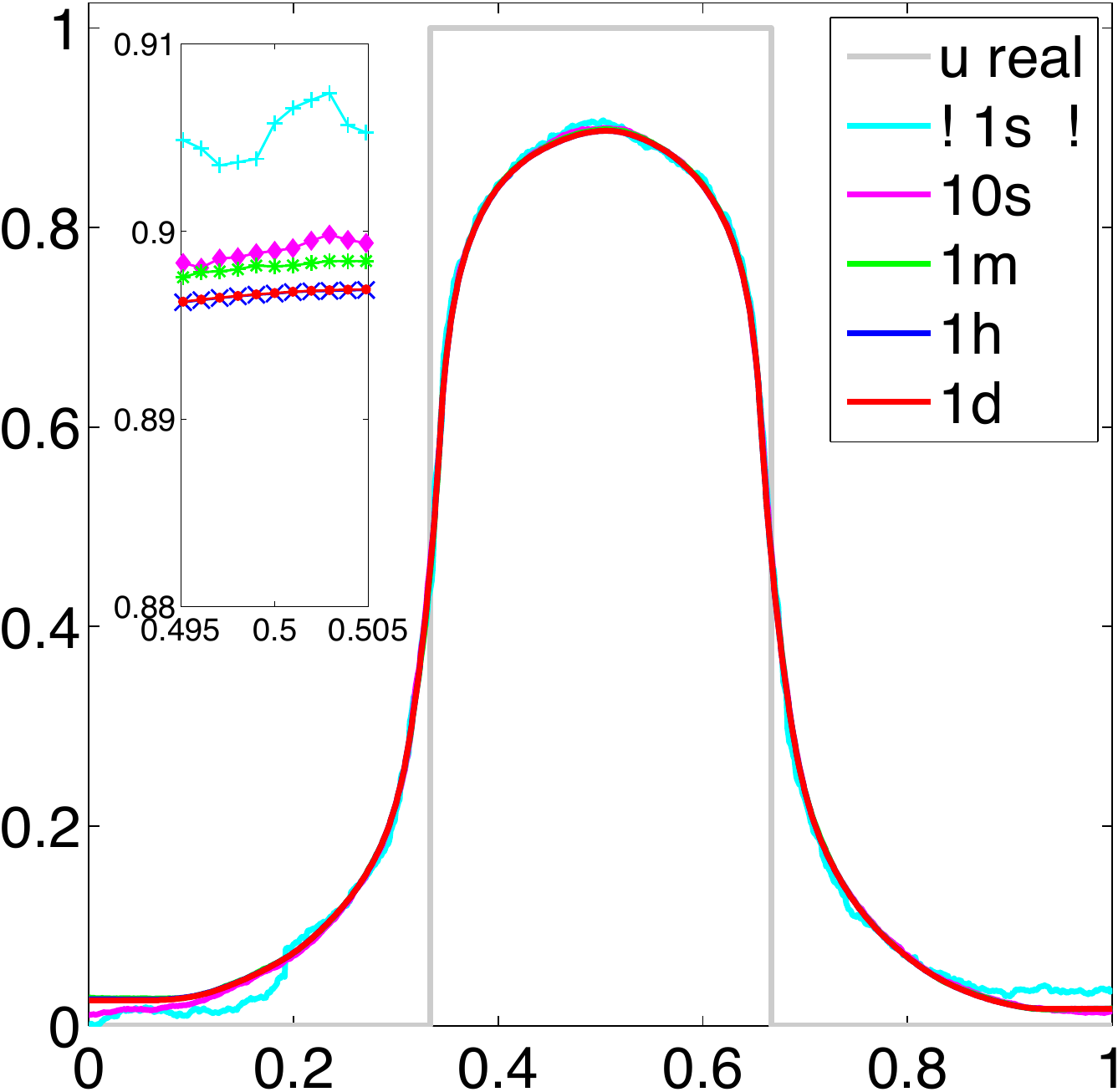}}
\subfigure[][RnGibbsO3\label{subfig:visu_RnGibbsO3}]{\includegraphics[width=0.3\textwidth]{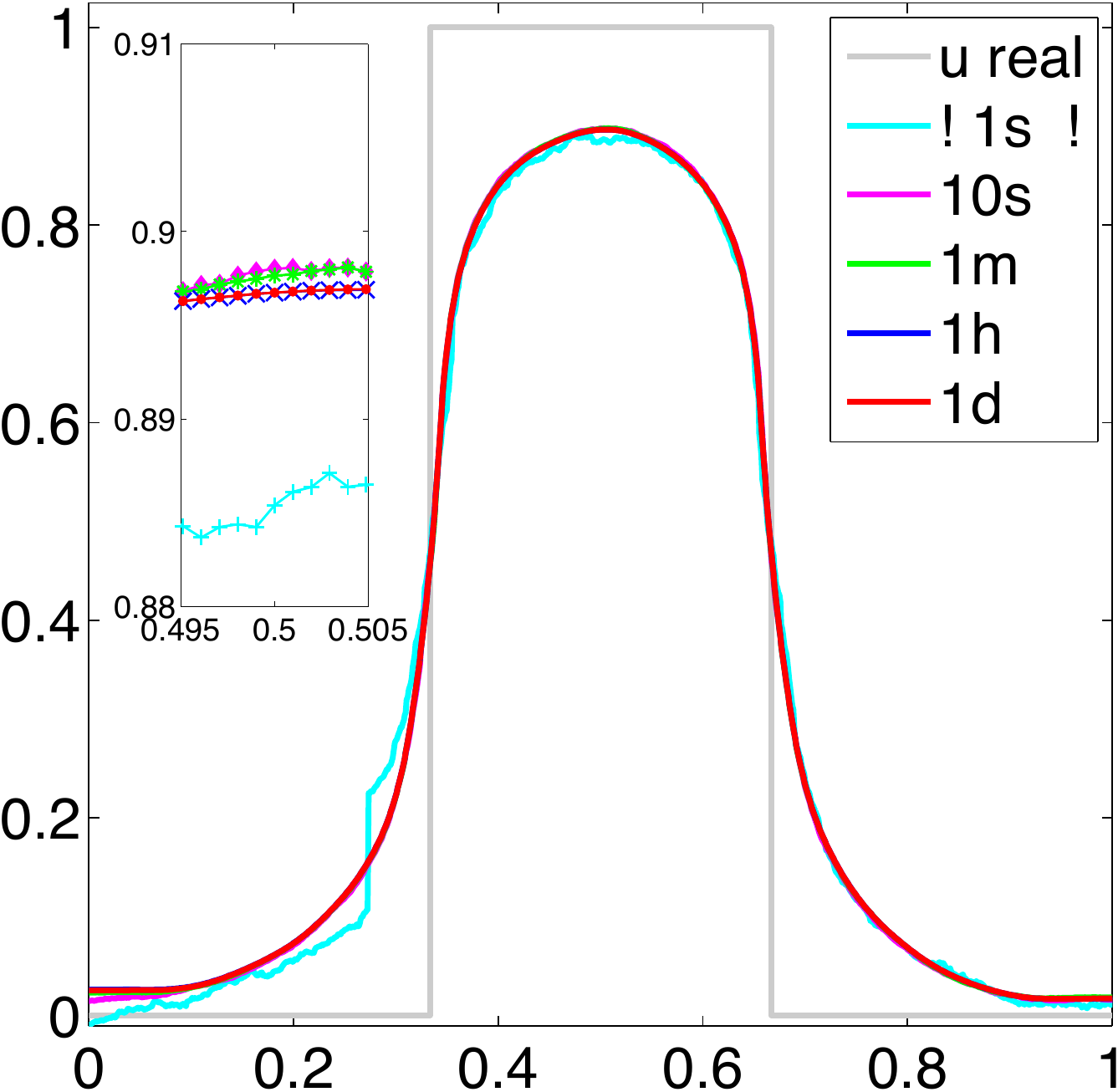}}
\subfigure[][RnGibbsO7\label{subfig:visu_RnGibbsO7}]{\includegraphics[width=0.3\textwidth]{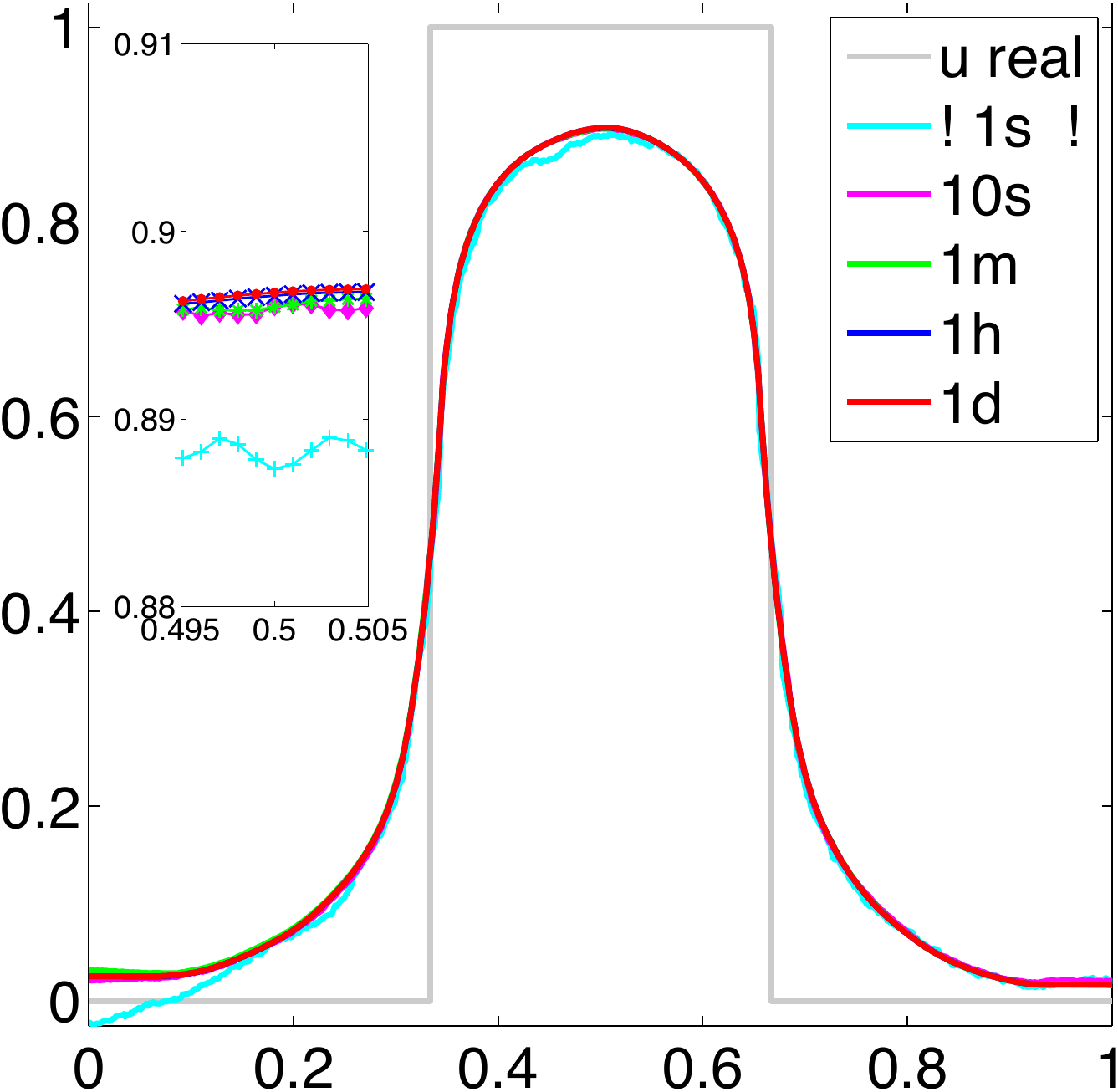}}\\
\subfigure[][SysGibbs\label{subfig:visu_SysGibbs}]{\includegraphics[width=0.3\textwidth]{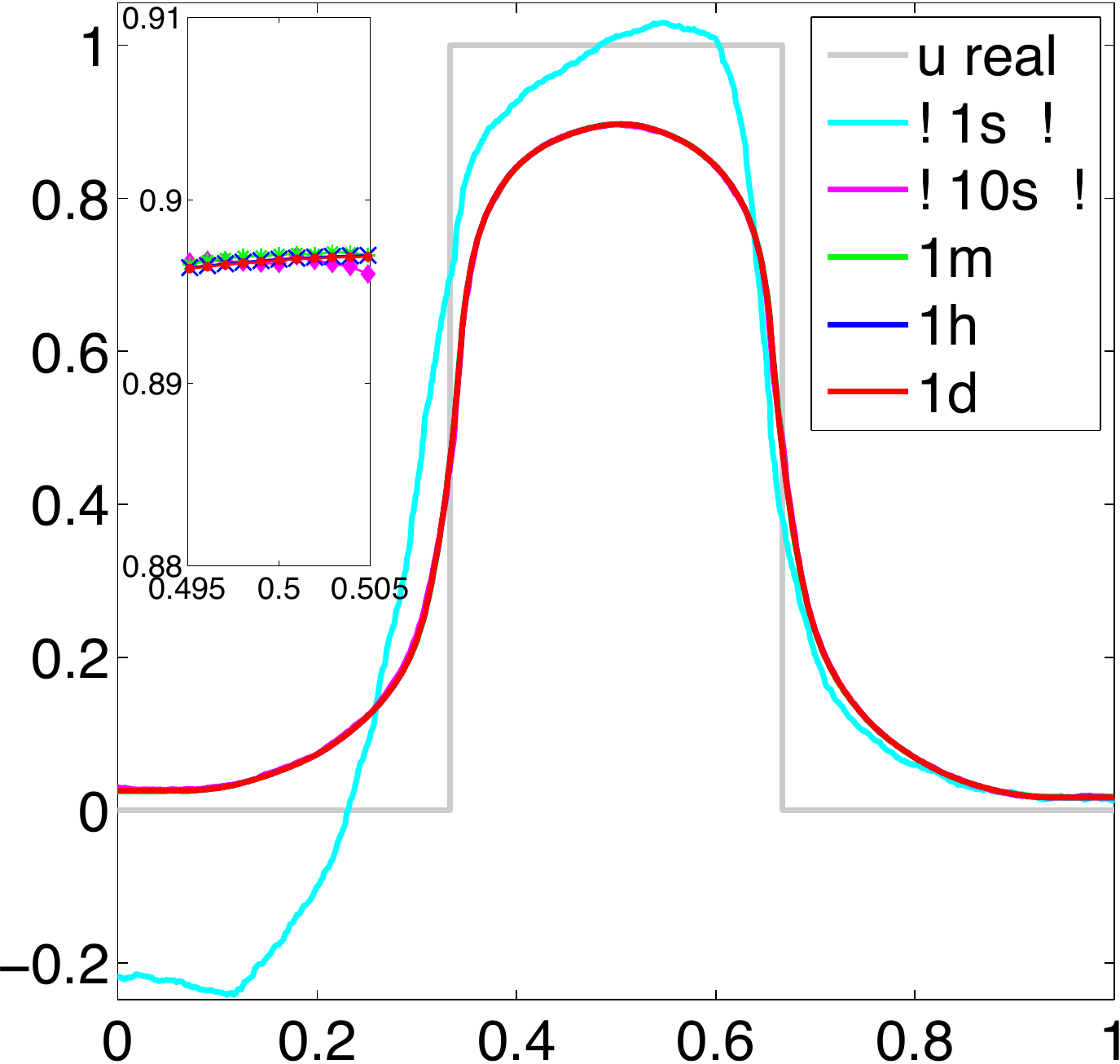}}
\subfigure[][SysGibbsO3\label{subfig:visu_SysGibbsO3}]{\includegraphics[width=0.3\textwidth]{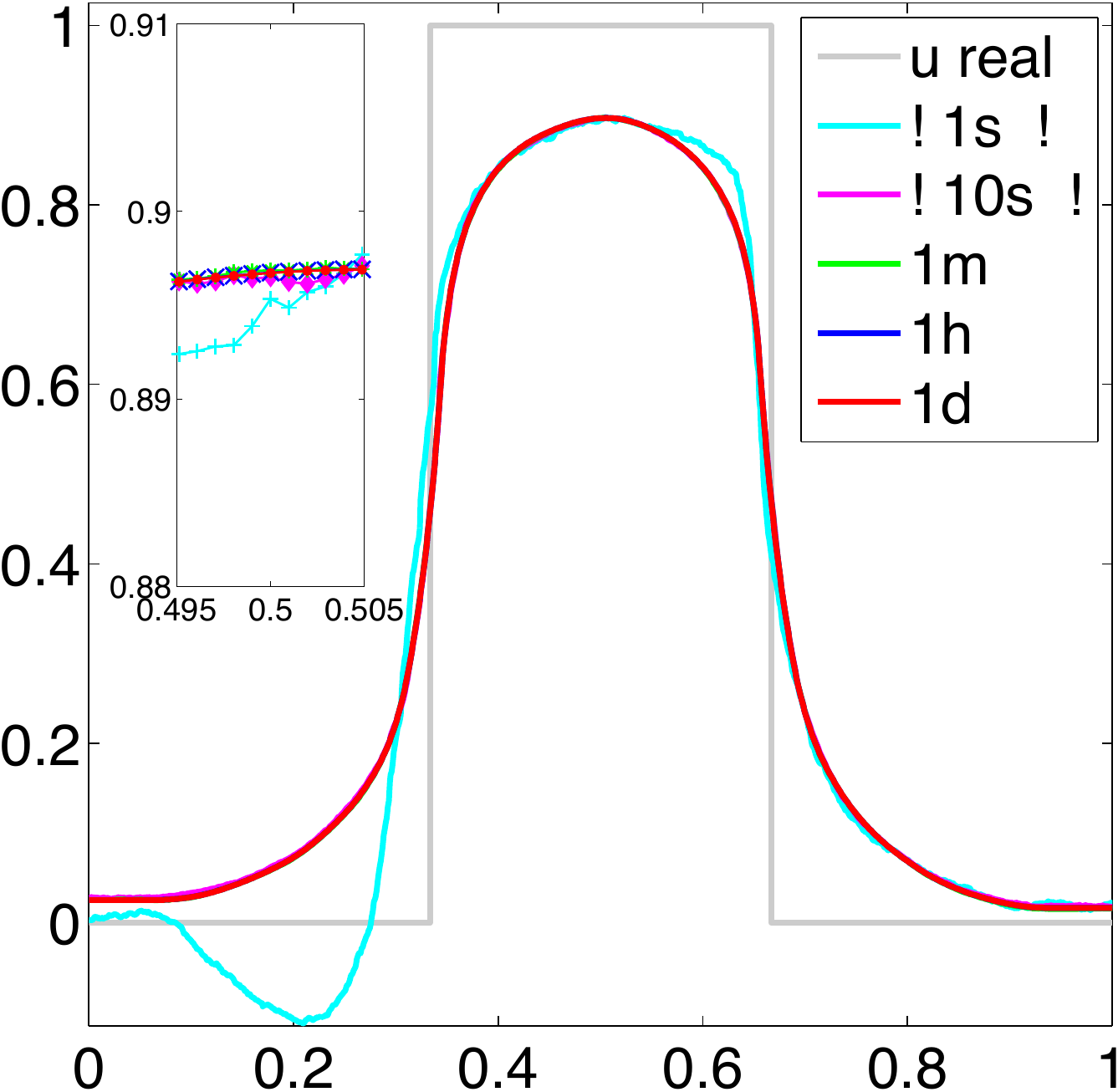}}
\subfigure[][SysGibbsO7\label{subfig:visu_SysGibbsO7}]{\includegraphics[width=0.3\textwidth]{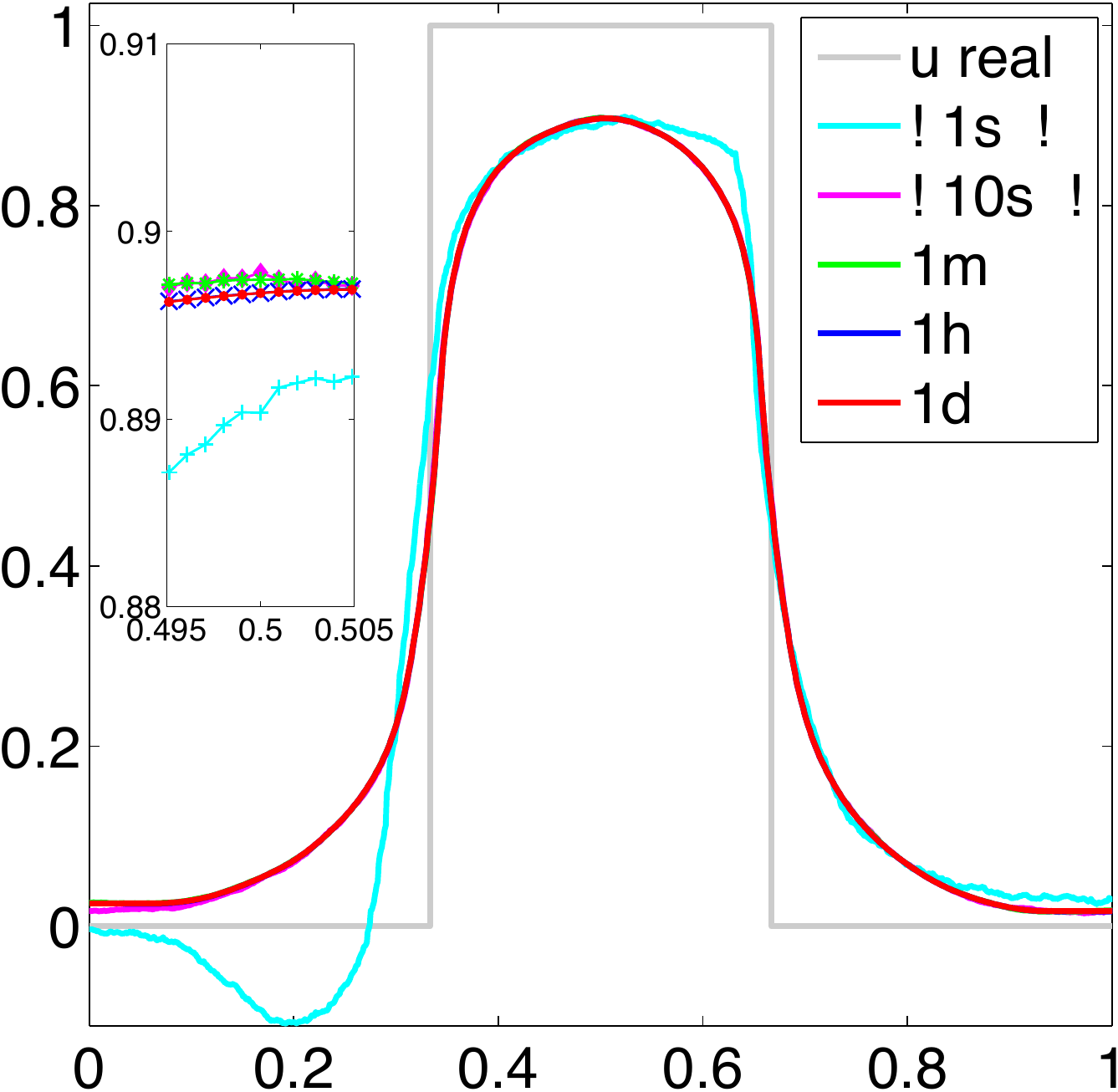}}
  \caption[]{Visual impressions of the CM estimate obtained by a sampler after a certain computation time ($n$ = 1023). Computation times, which are framed by exclamation marks, indicate that the burn-in time listed in Table \ref{tbl:burn-in} was not reached yet. For the Gibbs samplers, a zoom into the interval $[0.495,0.505]$ is added.\label{fig:visu}}
\end{figure}

\begin{figure}[hbt]
   \centering
\includegraphics[width = \textwidth]{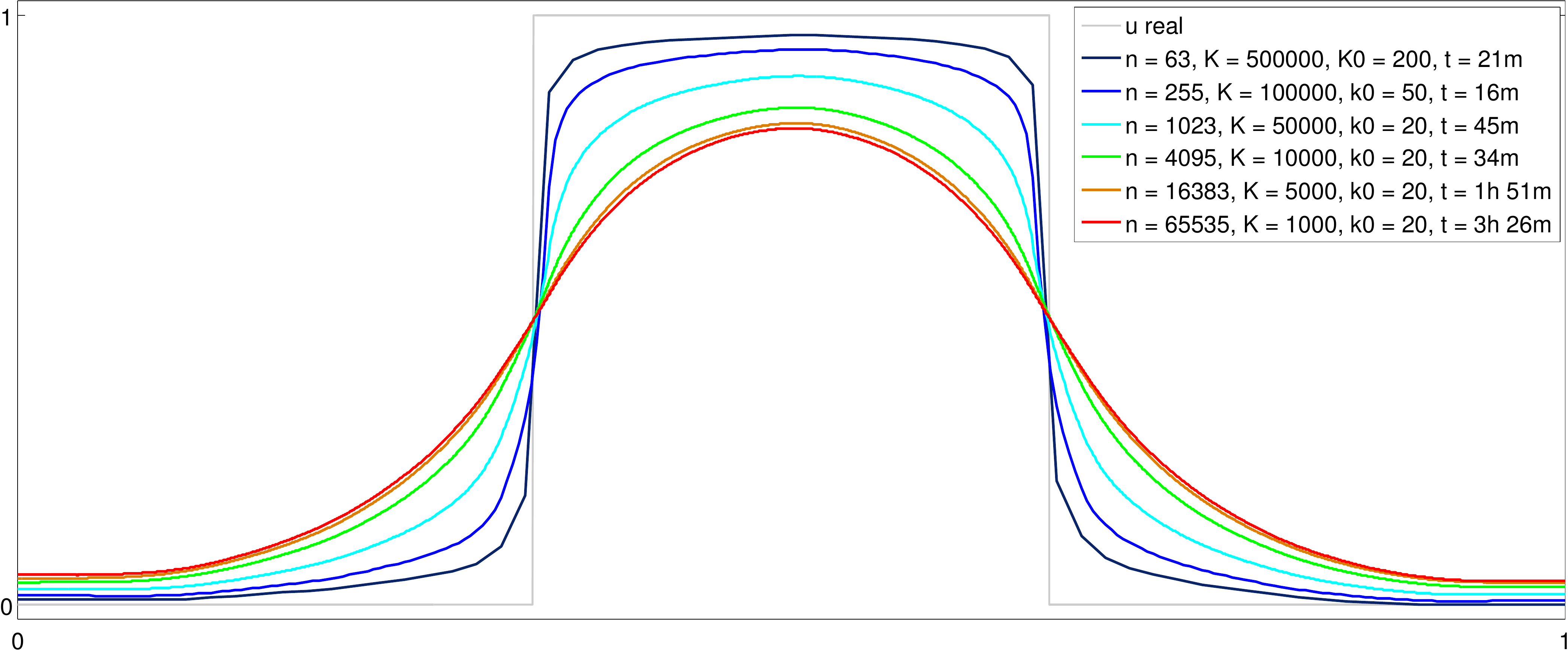}
\caption{CM estimates for growing $n$ and $\lambda_n \propto  \sqrt{n+1}$. Sampling was performed with the RnGibbs sampler.  \label{fig:PriorConv}}
\end{figure}

To get a visual impression of the sampling results, CM estimates are computed using the different samplers at different computation times for $n = 1023$, $\lambda_n = 25 \cdot \sqrt{n+1}$. The computation times examined are $t^* =$ 1 s, 10 s, 1 minute, 1 hour and  1 day, respectively. Practically, a long chain with $K_0 = 0$ was generated and sub-chains corresponding to all the samples drawn before $t^*$ were extracted. Then, CM estimates were computed from the sub-chains by discarding $K_0^* = \min(K_0, 0.5 \cdot K^*)$ burn-in samples, where $K_0$ are the burn-in steps listed in Table \ref{tbl:burn-in}, and $K^*$ denotes the number of samples in the subchain. The results are shown in Figure \ref{fig:visu}.  
We have to emphasize that we did not chose to show the CM estimate for the TV prior because the reconstruction is convincing. In fact, as explained in Section \ref{subsec:EdgPreBayInv} they are extremly smooth compared to the corresponding MAP estimates and, thus, bad reconstructions of the discontinuous $\tilde{u}$ . However, this smoothness is very useful for gaining a visual impression of the convergence and the properties of the different sampling schemes: The CM estimate computed from the chain converged once it is smooth. To demonstrate the capabilities of the newly developed Gibbs samplers for the practical use, we also examine the theoretical questions addressed in \cite{LaSi04}. For the choice of $\lambda_n \propto  \sqrt{n+1}$, the TV prior converges to a smoothness prior. To support this finding with numerical simulations the CM estimate was computed for $n = 63, 255, 1023, 4095$ in \cite{LaSi04} using the MH-Ncom sampler. Although the whole computation took about a month of time on a desktop PC equipped with a 2.8 GHz single core CPU, the authors admitted that the results were only partly satisfying. In Figure \ref{fig:PriorConv}, we show the CM estimate computed for $n = 63, 255, 1023, 4095, 16383, 65535$ using the RnGibbs sampler on a comparable CPU. Again, the CM estimates are only shown because the increasing smoothness of the CM estimates for growing $n$ allows for the visual inspection of the chain convergence.  

\subsubsection{Normal vs. ORR Gibbs Sampling} \label{subsubsec:NormVsOOR}
\begin{figure}[hbt]
   \centering
\subfigure[][RandGibbsO$\ast$ \label{subfig:RandOORAna}]{\includegraphics[height=5.4cm]{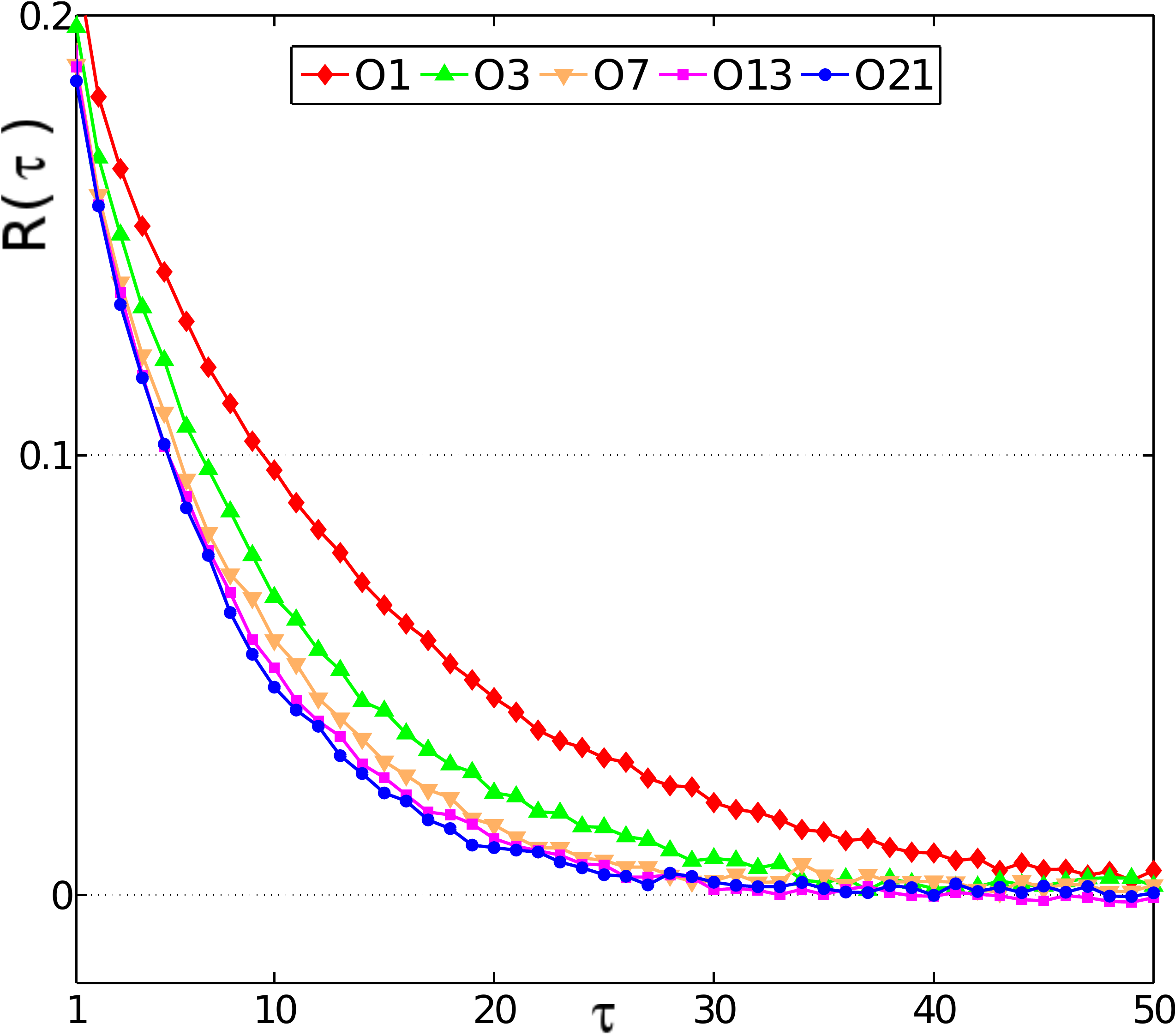}}
\subfigure[][SysGibbsO$\ast$ \label{subfig:SysOORAna}]{\includegraphics[height=5.4cm]{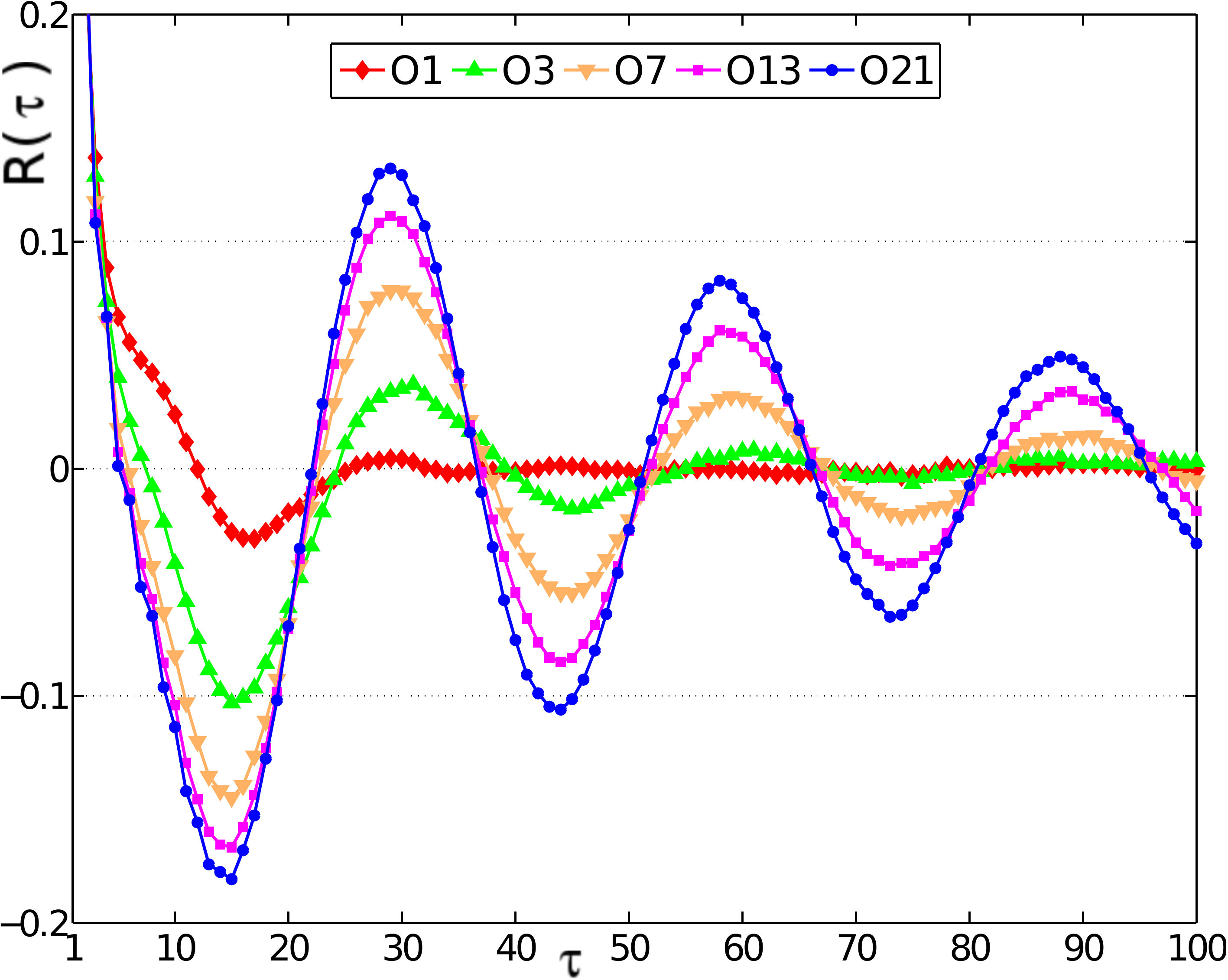}}
\caption{Autocorrelation plots for random and systematic scan Gibbs sampling using different values for oriented overrelaxation parameter $N_O$. Red \mbox{$\blacklozenge$} : O1 (i.e., no overrelaxation at all); green \mbox{$\blacktriangle$} : O3; orange \mbox{$\blacktriangledown$} : O7; pink \mbox{$\blacksquare$} : O13; blue \mbox{$\bullet$} : O21.
   \label{fig:Rand_vs_Sys}}
\end{figure}

Using oriented overrelaxation removes a certain amount of randomness from the generated chains. This may lead to a faster exploration of the posterior distribution, but can also enhance non-ergodic tendencies of the sampling approaches. In Figure \ref{fig:Rand_vs_Sys}, more detailed autocorrelation plots for the 1D scenario using $n = 1\,023$, $\lambda = 800$ are shown. Whether oriented overrelaxation is practically advantageous relies on the computational cost of sampling from the single component density \eref{StdDens} costs compared to the other computation steps. Using oriented overrelaxation to sample from \eref{StdDens} takes roughly twice as much computation time as not using it, almost independent from $N_O$ (when using Algorithm \ref{algo:ImplOOR}). If the other computation steps in the whole sampling scheme take way more time (i.e., the computation of $b$, see Section \ref{subsec:ImpGibbsL1}), this extra computational cost is negligible. In our scenario, the ratio between the total computation time per sample $t_s$ for the RnGibbsO7 and the RnGibbs sampler varies considerably. It drops from 1.93 for the fast but memory consuming implementation (see Section \ref{subsec:ImpGibbsL1} ) to compute $b$ and $n=63$ to 1.02 for the slower implementation to compute $b$ and $n = 65\,535$.

\subsection{Image Deblurring with Impulse Prior in 2D} \label{subsec:ImDe2D}
\subsubsection{Setting} \label{subsubsec:Set2D}
\begin{figure}[hbt]
   \centering
\subfigure[][  \label{subfig:cmap2D}]{\includegraphics[height = 0.45 \textwidth]{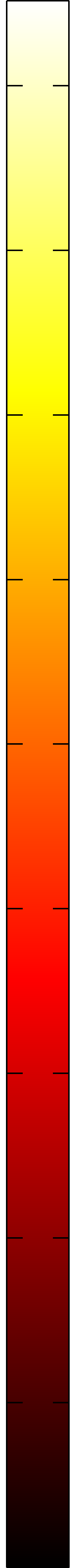}}
\subfigure[][ \label{subfig:Real2}]{\includegraphics[height = 0.45\textwidth]{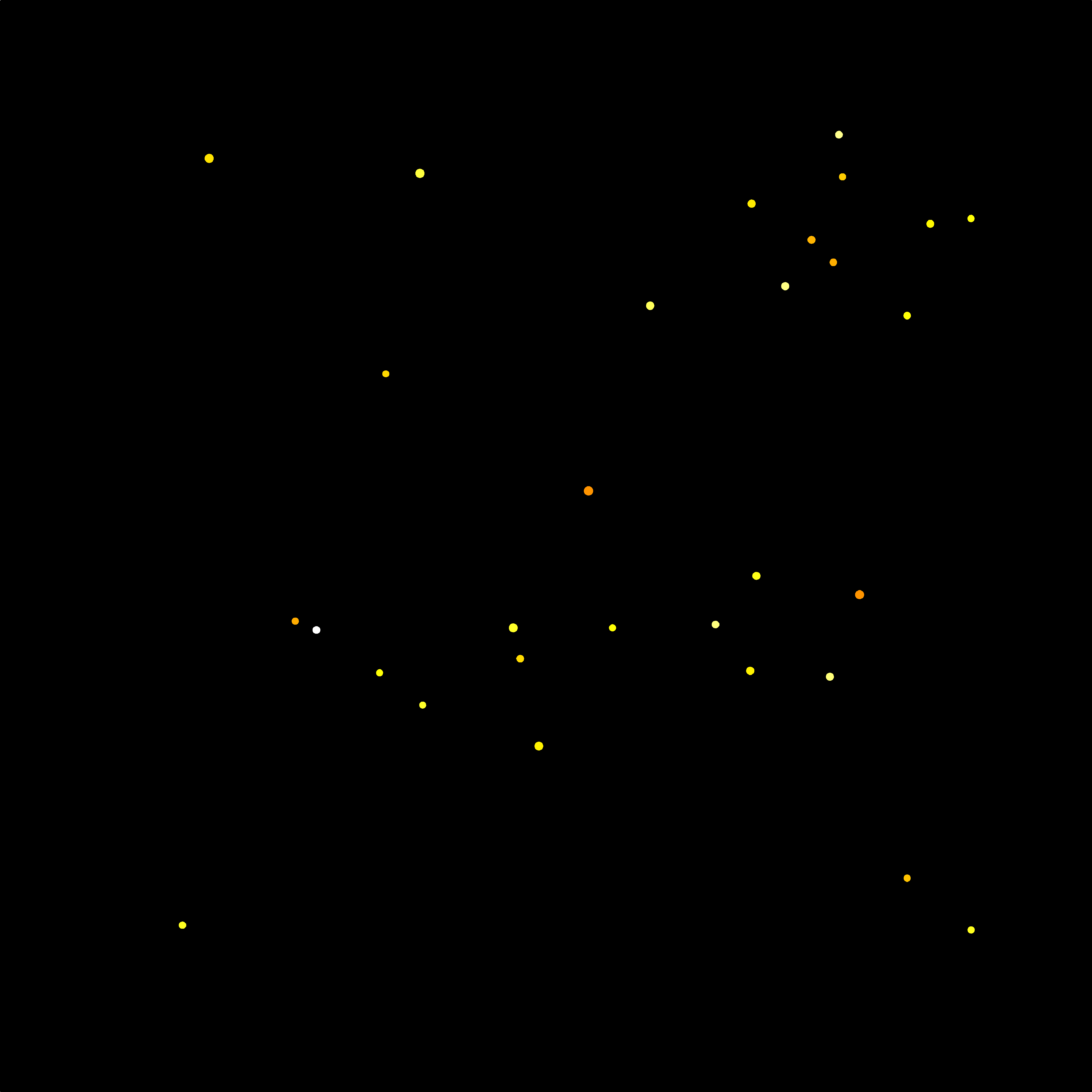}}
\subfigure[][  \label{subfig:Data2}]{\includegraphics[height = 0.45 \textwidth]{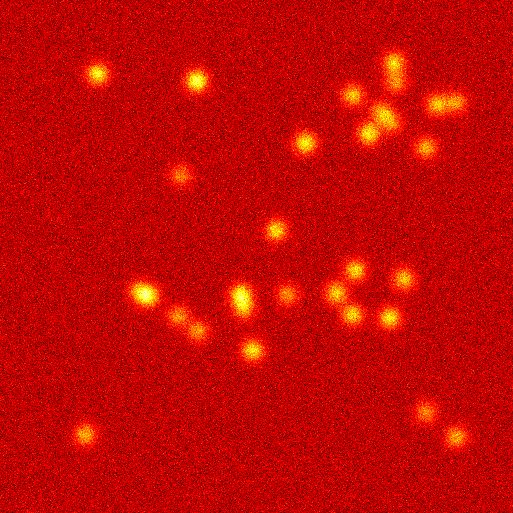}}
\caption{Left: The colour scale used in all 2D images. Middle: The unknown function $\tilde{u}$. Right: The measurement data $m$. }
   \label{fig:RealData2}
\end{figure}

As a second example, we consider 2D image deblurring with a simple L1 prior,  i.e., $D = I_n$  (also called  \emph{impulse prior}). The unknown intensity function $\tilde{u}:[0,1]\times [0,1] \rightarrow \mathbb{R}$ is shown in Figure \ref{subfig:Real2}. It consists of a couple of circular spots of constant intensity whose radii and intensities slightly vary between single spots. The forward mapping is given by a convolution with a Gaussian kernel with standard deviation of 0.015. Measurement data is generated by integrating the resulting convoluted image over $513 \times 513$ regular pixel and adding noise. The relative noise level is 0.1, i.e., the standard deviation $\sigma$ of the measurement noise is 0.1 times the maximal intensity of the noiseless signal. The resulting measurement data is shown in Figure \ref{subfig:Data2}. The image will be reconstructed on the same pixel grid used for the measurement using Neumann boundary conditions, thus, the dimension of the unknowns $n$ is $511^2 = 261\,121$. To avoid an inverse crime, the grid used for the generation of the measurement data was 4 times finer.
For the MH samplers, the same on-line adaptation of $\kappa$ was used as explained in Section \ref{subsubsec:Pre}. The update intervals and the up- and down-scaling factors have been chosen carefully to optimize the performance of the samplers while keeping the adaptation stable, i.e., monotonic. The forward mapping is implemented using ffts. For the Gibbs sampler, we note that $V = I_n$  and that \eref{eq:2DimplB} simplifies to 
\begin{equation}
 (\psi_i^t \, \Psi_{[-i]}) \, \xi_{[-i]} =  \frac{1}{2 \sigma^2}  \left[ \left(A^t A \right) \cdot \xi \right]_i  - \xi_i  \|  \psi_i \|^{2}_{2},
\end{equation}
which can be implemented in a efficient, direct way. 

\subsubsection{Visual Results} \label{subsubsec:VisRes2D}
\begin{figure}[hbt]
   \centering
\subfigure[][MH-Iso, 1h \label{subfig:MH-Iso_1h}]{\includegraphics[width = 0.328\textwidth]{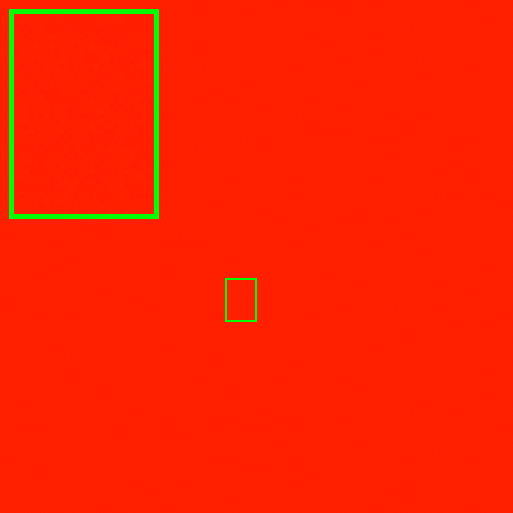}}
\subfigure[][MH-Iso, 5h \label{subfig:MH-Iso_5h}]{\includegraphics[width = 0.328 \textwidth]{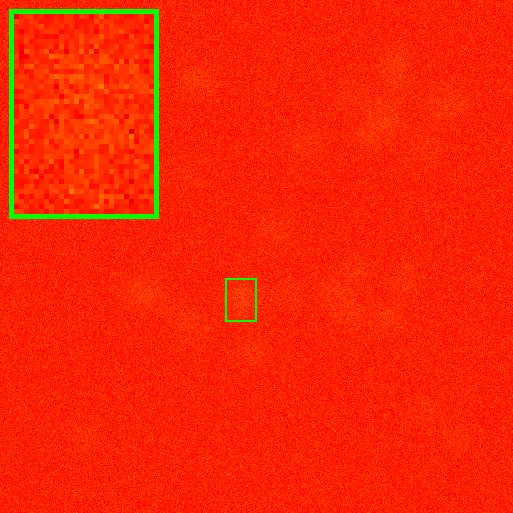}}
\subfigure[][MH-Iso, 20h \label{subfig:MH-Iso_20h}]{\includegraphics[width = 0.328 \textwidth]{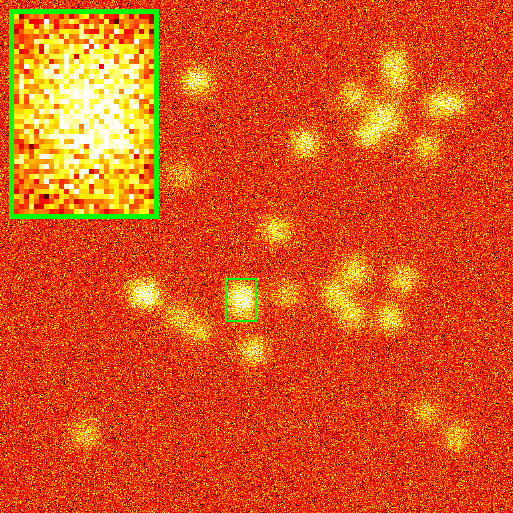}}\\
\subfigure[][MH-Si, 1h \label{subfig:MH-Si_1h}]{\includegraphics[width = 0.328\textwidth]{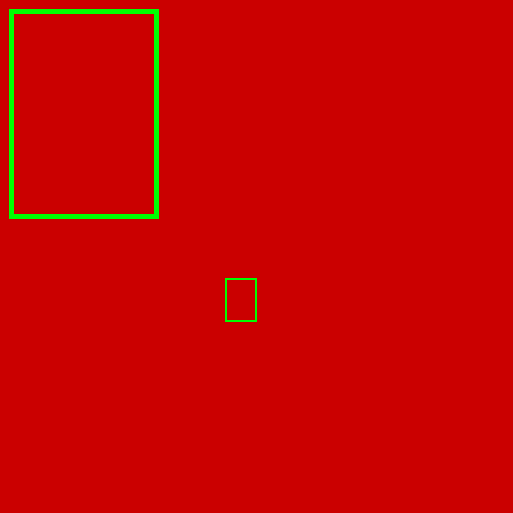}}
\subfigure[][MH-Si, 5h \label{subfig:MH-Si_5h}]{\includegraphics[width = 0.328 \textwidth]{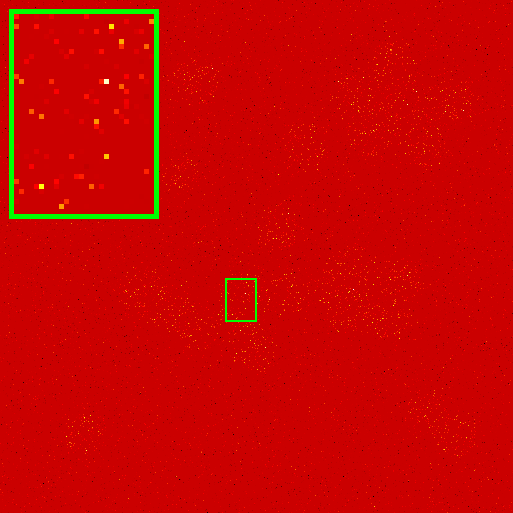}}
\subfigure[][MH-Si, 20h \label{subfig:MH-Si_20h}]{\includegraphics[width = 0.328 \textwidth]{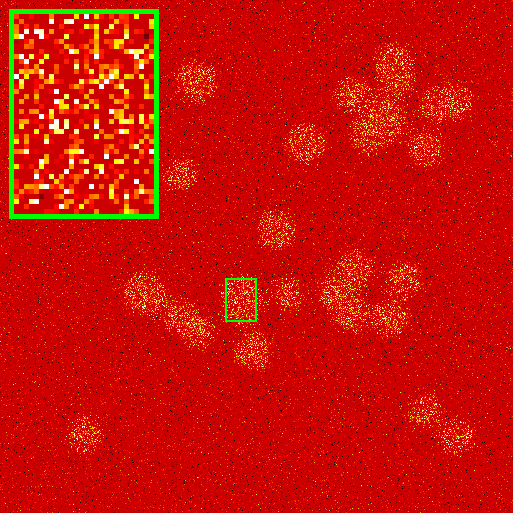}}\\
\caption{Visual results for image deblurring with an impulse prior in 2D, part 1. In the upper left corner of each subfigure, a zoom into the marked area in the original figure is shown. The scaling used in these images is explained in the text.}
   \label{fig:Visu2Da}
\end{figure}

\begin{figure}[hbt]
   \centering
\subfigure[][RnGibbs, 1h \label{subfig:RnGibbs_1h}]{\includegraphics[width = 0.328\textwidth]{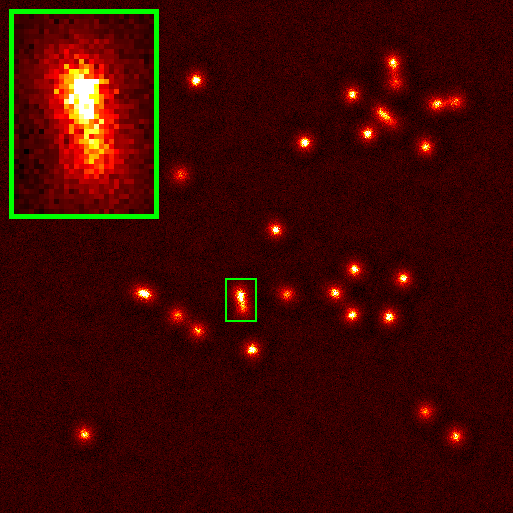}}
\subfigure[][RnGibbs, 5h \label{subfig:RnGibbs_5h}]{\includegraphics[width = 0.328 \textwidth]{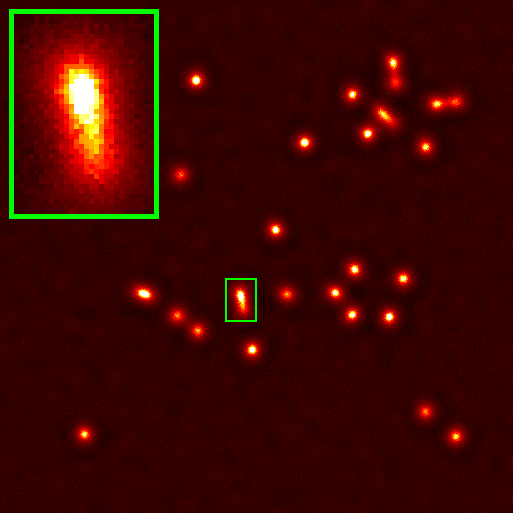}}
\subfigure[][RnGibbs, 20h \label{subfig:RnGibbs_20h}]{\includegraphics[width = 0.328 \textwidth]{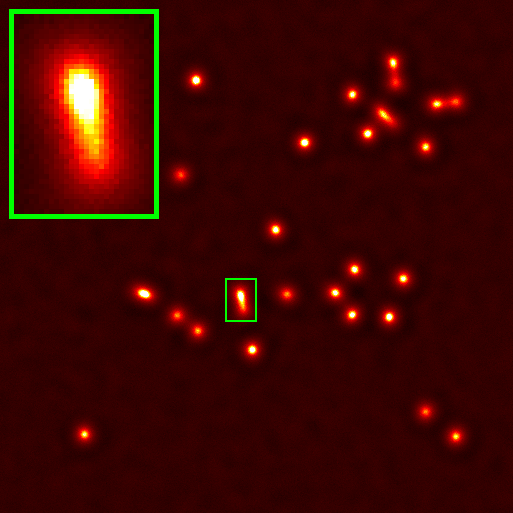}}\\
\subfigure[][RnGibbsO7, 1h \label{subfig:RnGibbsO7_1h}]{\includegraphics[width = 0.328\textwidth]{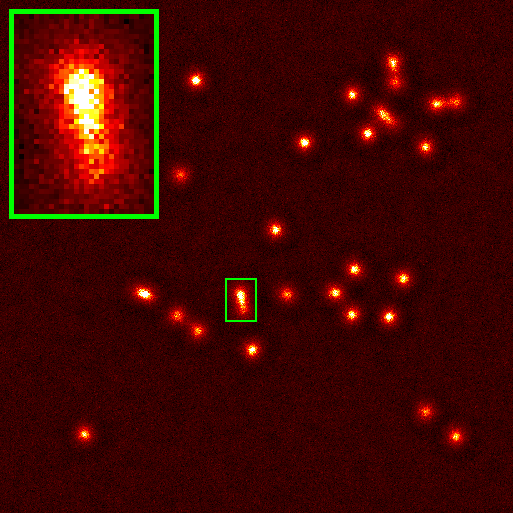}}
\subfigure[][RnGibbsO7, 5h \label{subfig:RnGibbsO7_5h}]{\includegraphics[width = 0.328 \textwidth]{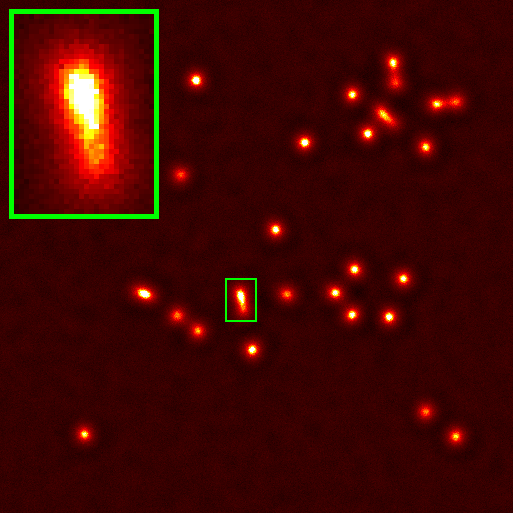}}
\subfigure[][RnGibbsO7, 20h \label{subfig:RnGibbsO7_20h}]{\includegraphics[width = 0.328 \textwidth]{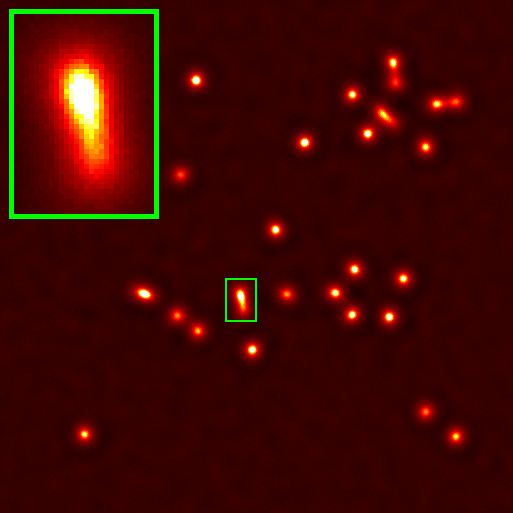}}\\
\caption{Visual results for image deblurring with an impulse prior in 2D, part 2. In the upper left corner of each subfigure, a zoom into the marked area in the original figure is shown. The scaling used in these images is explained in the text.}
   \label{fig:Visu2Db}
\end{figure}

\begin{figure}[hbt]
   \centering
\subfigure[][SysGibbs, 1h \label{subfig:SysGibbs_1h}]{\includegraphics[width = 0.328\textwidth]{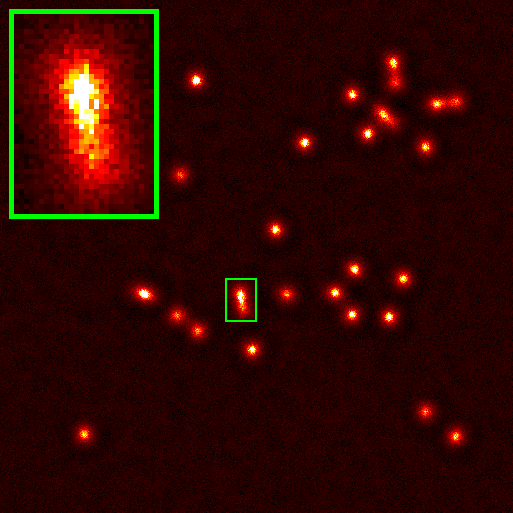}}
\subfigure[][SysGibbs, 5h \label{subfig:SysGibbs_5h}]{\includegraphics[width = 0.328 \textwidth]{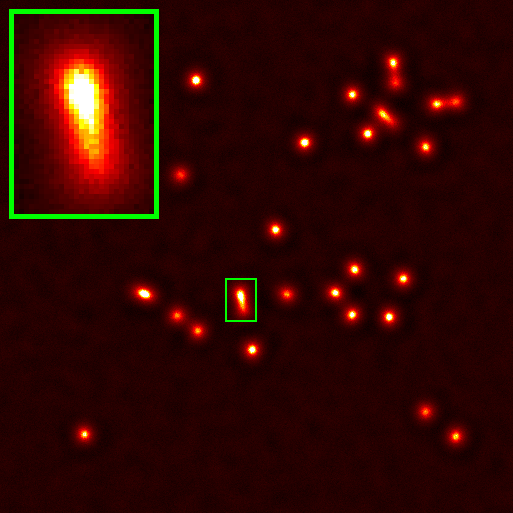}}
\subfigure[][SysGibbs, 20h \label{subfig:SysGibbs_20h}]{\includegraphics[width = 0.328 \textwidth]{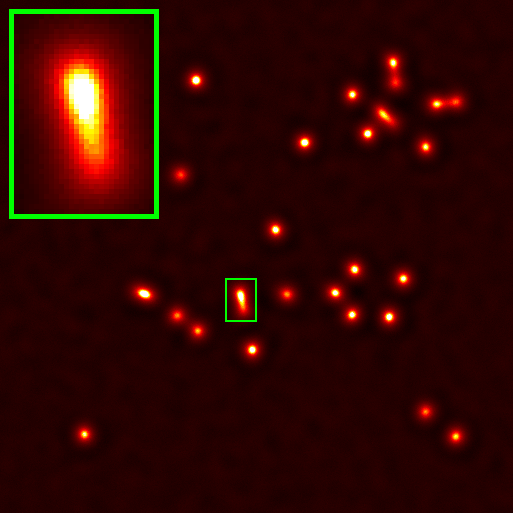}}\\
\subfigure[][SysGibbsO7, 1h \label{subfig:SysGibbsO7_1h}]{\includegraphics[width = 0.328\textwidth]{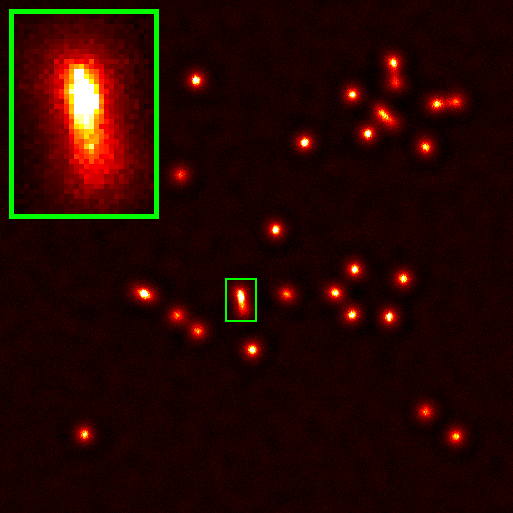}}
\subfigure[][SysGibbsO7, 5h \label{subfig:SysGibbsO7_5h}]{\includegraphics[width = 0.328 \textwidth]{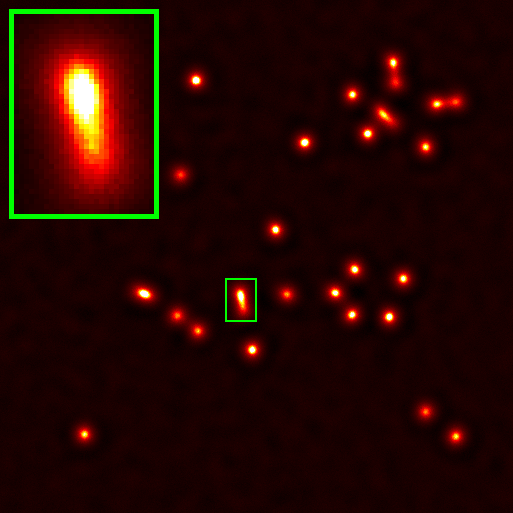}}
\subfigure[][SysGibbsO7, 20h \label{subfig:SysGibbsO7_20h}]{\includegraphics[width = 0.328 \textwidth]{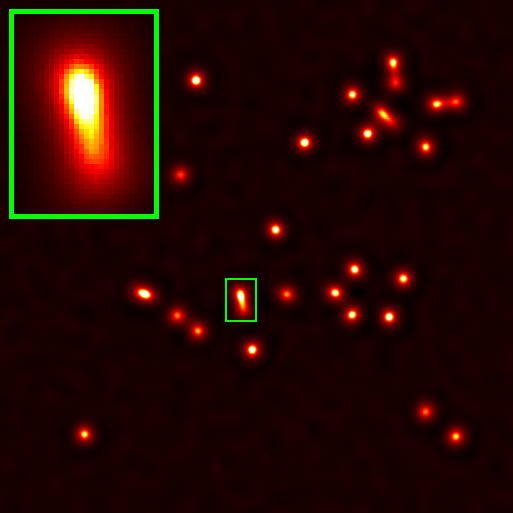}}\\
\caption{Visual results for image deblurring with an impulse prior in 2D, part 3. In the upper left corner of each subfigure, a zoom into the marked area in the original figure is shown. The scaling used in these images is explained in the text.}
   \label{fig:Visu2Dc}
\end{figure}

The practical procedure to compute visual  results is identical to the one used in Section \ref{subsubsec:VisRes1D}. In Figures \ref{fig:Visu2Da}-\ref{fig:Visu2Dc}, the CM estimates computed after 1, 5 and 20 hours are shown. The results of MH-Ncomp, RnGibbsO3 and SysGibbsO3 are omitted here. Choosing a good scaling to compare the results for a single sampling method is not easy because of outliers in the 1h image. These outliers would either lower the contrast if a simple linear min-max scaling based on all images is chosen or would lead to the impression that constant regions are growing if an individual scaling for each image is used. The scaling we used is generated in the following way: For each method, we merged and sorted the pixel values of all three CM estimates. From this sorted set, the smallest and largest values are discarded, using 0.1\% and 99.9\% as thresholds, respectively. A linear min-max scaling is generated from the remaining values, and the discarded values are mapped to the beginning and end to this color scale, respectively.\\ 
An examination of $\log[p(u_i|m)]$ was again used to determine the burn-in steps $K_0$. For the Gibbs samplers, the burn-in steps were between 12 and 30 and in all images shown, the burn-in phase was already completed. For the MH samplers, the plots of $\log[p(u_i|m)]$ suggested that even after 20 hours of computation, the chain was still far away from the central parts of the posterior, which is also evident from the CM estimates.

\section{Discussion and Conclusions}\label{sec:Discussion}
\subsection{MH-Samplers} \label{subsec:DisMH-Samplers}
For the  specific scenario we examined, the efficiency of the basic MH samplers decreases when either the influence of the L1-type prior increases (i.e., $\lambda$ is increased) or the dimension of the unknowns, $n$, is increased. Figures \ref{fig:acf_n63_MH}, \ref{fig:tacf_n63_MHvsGibbs}, \ref{fig:acf_nVar_MHvsGibbs}  and \ref{fig:tacf_nVar_MHvsGibbs} and Tables \ref{tbl:burn-in} and \ref{tbl:lag1} clearly document this. The largest number of unknowns examined was $n = 1\,023$, which is still moderate for typical inverse problem scenarios. However, for this number of unknowns, both the burn-in time and the time to decrease the autocorrelation of a new sample below 1\% are in the order of a few hours. In Figures \ref{subfig:visu_MH-Iso} - \ref{subfig:visu_MH-Si}, this is visualized by the slow convergence of computed CM estimates in the 1D scenario. A (visually) satisfactory result is only obtained after 1 day of computation time. In the 2D example, no satisfactory result could be obtained, even after 20 hours of computation time (see Figure \ref{fig:Visu2Da}). The examination of $\log[p(u_i|m)]$ suggested that the computation required to obtain such a result is of orders larger. In total, our detailed studies support the empirical findings of former applications of basic MH-samplers to L1-type priors, see, e.g., \cite{LaSi04,KoLaNiSi12}. 

\subsection{Gibbs-Samplers} \label{subsec:DisGibbs-Samplers}
We again stress that the Gibbs samplers we proposed and examined have to be considered as very basic variants of Gibbs sampling (cf. Section \ref{subsec:PosInfMCMSam}). This makes it even more surprising that for these samplers show totally different trends compared to the MH-samplers. 
\paragraph{Random Scan Gibbs Samplers:} For RnGibbs, RnGibbsO3 and RnGibbsO7, the required burn-in steps stay constant when increasing $\lambda$ and clearly decrease when increasing $n$, cf. Table \ref{tbl:burn-in}. Even as the computational costs of drawing a new sample increases with $n$, this effect keeps the computational time to draw the required number of burn-in steps almost constant. Figures \ref{fig:acf_n63_Gibbs}, \ref{fig:tacf_n63_Gibbs} and Table \ref{tbl:lag1} show that the decay of $R(\tau)$ and $R^*(t)$ is even faster for increasing $\lambda$. From Figures \ref{fig:acf_nVar_MHvsGibbs}, \ref{fig:acf_nVar_Gibbs} and Table \ref{tbl:lag1} we see that for increasing $n$, this is also true for $R(\tau)$. For $R^*(t)$, we see in Figure \ref{fig:tacf_nVar_MHvsGibbs} that for large $n$ the temporal decay cannot further decrease. This is a normal saturation effect because the the autocorrelation decrease is bounded. It would even occur for an i.i.d. series of $n$ dimensional random variables, if the computation time would increase with $n$. The visual results (see Figures \ref{subfig:visu_RnGibbs} - \ref{subfig:visu_RnGibbsO7}, \ref{fig:Visu2Db},\ref{fig:Visu2Dc}) clearly support these findings. Especially the short burn-in times are noticeable. In both scenarios, the CM estimate using the shortest computation time already represents the most important features of the final solution. Using oriented overrelaxation in combination with random scan Gibbs sampling does not seem to lead to any problems concerning the ergodicity of the chain. The autocorrelation plots in Figures  \ref{fig:acf_n63_Gibbs}, \ref{fig:acf_nVar_Gibbs} \ref{subfig:RandOORAna} are still monotonic and positive. Clearly, the decay of $R(\tau)$ is faster using overrelaxation. Concerning $R^*(t)$, overrelaxation is only effective, if the additional computational cost is negligible compared to other parts of the sampling process. 
\paragraph{Systematic Scan Gibbs Samplers:} For SysGibbs, SysGibbsO3 and SysGibbsO7, the results are less clear. At first glance, the trends in their results seem to be rather similar to the random scan samplers and within a direct comparison, they often seem to outperform them, see, e.g., Table \ref{tbl:lag1}. However, Figures \ref{fig:acf_n63_Gibbs}, \ref{fig:acf_nVar_Gibbs} and \ref{subfig:SysOORAna} show that for growing $n$ and $N_O$, the plots of $R(\tau)$ start to oscillate and are clearly negative in some areas. It seems that in combination with the TV prior, the subsequent update of neighboring increments leads to non-ergodic tendencies in the sampling procedure. These tendencies are amplified when using oriented overrelaxation. Producing anti-correlated samples may in fact advantageous for certain tasks \cite{Li08}. However, we would, in general, not advise to use systematic Gibbs sampling. The additional computational cost of drawing a random component to update in each step is small. In contrast, a sampler that relies on a non-ergodic mapping may produce unpredictable results for certain tasks.

\subsection{General} \label{subsec:GenDis}
There are multiple reasons for the loss of performance of the basic MH samplers compared to the basic Gibbs samplers in the specific scenarios we examined. The crucial part for an MH sampler is the design of a good proposal distribution. As explained in Section \ref{subsec:PosInfMCMSam}, the basic MH samplers we applied are ``black-box sampler'' algorithms. In the design of their proposal distributions, no specific information about the posterior was taken into account. In return, they usually exhibit very fast computation times. The standard proposal distributions we used are designed to sample from low dimensional, Gaussian-like distributions. However, high dimensional posteriors from sparsity promoting priors have very different properties. Standard MH-samplers have to take very small steps to obtain a good acceptance rate. This leads to long burn-in times and a slow decrease in autocorrelation. The situation is similar with optimization algorithms used for MAP estimation (cf., Section \ref{subsec:GenSetBayFor}). Black-box optimization algorithms that only rely on evaluating the objective function (i.e, $\log[p(u|m)]$)  are usually too slow when applied to specific, high dimensional posteriors. The basic Gibbs samplers we proposed incorporate more posterior-specific information into the sampling procedure at the costs of a larger computation time. The conditional single component densities, which can be regarded as optimal transition kernels, are computed and sampled from explicitly. This small extra amount of incorporating problem specific information already seems to be sufficient to generate very promising sampling procedures for high dimensional Bayesian inversion using L1-type priors (the dimensions of the unknowns used in Figure \ref{fig:PriorConv} and the 2D scenario are far beyond any previously reported use of MCMC for L1-type inverse problems). In general, both sampling techniques have advantages and disadvantages and will outperform the other given a specific scenario.\\
For the very reason that we only used very basic Gibbs samplers, our results also challenge common beliefs about the feasibility of MCMC sampling in high dimensional inverse problems in general. We showed that MCMC schemes are not in general slow and scale bad with increasing dimension. We rather think that MCMC schemes for inverse problems are far less elaborate compared to optimization schemes up to now. \\
With regard to the corresponding optimization algorithms for MAP estimates, one possible reason for the superior performance of the single component Gibbs samplers might be the transformation of the posterior into the basis ${v_1,\ldots,v_n}$ (cf., \ref{subsec:ImpGibbsL1}). In this basis, the prior diagonalizes. It can be shown that the MAP estimate is sparse in this basis, i.e., many basis coefficients are exactly zero. Many optimization algorithms to compute the MAP estimate take advantage of this and perform better when transformed into that basis. One could argue that this might be the case for the sampling procedures as well, and that a fair comparison between MH and Gibbs samplers would need to transform the MH samplers into the basis ${v_1,\ldots,v_n}$ as well. However, there are reasons why this argument is not valid. The striking advantage of MH samplers is their simple, ``black-box''-like implementation. In practice, they are normally implemented in the most direct way and we stuck to that paradigm. In addition, in the 2D case, $v_i = e_i$, i.e., both MH and Gibbs samplers are already formulated in the right basis. However, this does not affect the bad performance of the MH samplers compared to the Gibbs samplers. But most importantly, while the MAP estimate is sparse, the CM estimate is not, and single samples from the posterior are not sparse as well. Theoretically, it has been shown in \cite{Lo08} that for denoising using a TV prior, the CM estimate is, in fact, never sparse. One can see this, e.g., in Figure \ref{fig:visu}. The estimates are neither sparse in the normal basis, nor in the increment basis. In Figures \ref{fig:Visu2Da}- \ref{fig:Visu2Dc}, one can clearly see that this is similar for the normal L1 prior.

\subsection{Outlook and Extensions} \label{subsec:Out}
In this first study, we only compared very basic variants of MH and Gibbs sampling. In the future, a comparison to more sophisticated variants of MH schemes such as \emph{delayed rejection} \cite{Mi01,HaLaMiSa06}, \emph{adaptive Metropolis} schemes \cite{HaSaTa05,HaLaMiSa06} or the \emph{t-walk} \cite{ChFo10} has to be undertaken. The most promising technique for our scenario might be to combine a tailored variant of delayed rejection with SCAM \cite{HaSaTa05}. Many of the sophisticated MH variants have been developed for the study of non-linear, computationally extensive and high dimensional inverse problems (see \cite{CuFoSu11} for a recent overview), i.e., situations where MH sampling is the only MCMC technique that can be applied. Improving the basic Gibbs sampling schemes used here by adding adaptive elements or optimal directions is far less developed until now \cite{ChFoPeSa12} and is an interesting future topic of research. A comparison between sophisticated (possibly adaptive, i.e., non-markovian) variants of MH and Gibbs sampling will also need a concrete application scenario since a general comparison by the measures used in this article is less meaningful.\\
We assumed that the variance $\sigma^2$ of the noise term (cf. Section \ref{subsec:GenSetBayFor}) is known exactly (or a good estimate is available), which is not always the case in practical applications. The Bayesian framework can account for the uncertainty of this model parameter as well: $\sigma^2$ is treated like the other unknowns $u$ (but assumed to be independent from them) and the available information about its typical values are expressed by a prior $p_{pr}(\sigma^2)$. An advantageous choice for $p_{pr}(\sigma^2)$ is given by the \emph{conjugate prior} w.r.t. $p_{li}(m|u,\sigma^2)$, which is the \emph{inverse gamma distribution} \cite{Ge06}:
\begin{equation}
p(x| \alpha, \beta) = \frac{\beta^{\alpha}}{\Gamma(\alpha)}\: x^{- \alpha - 1}  \exp\left( - \frac{\beta}{x}\right), \label{eq:InvGamma}\\
\end{equation}
with \emph{shape} and \emph{scale} parameters $\alpha$ and $\beta$ ($\Gamma$ denotes the gamma function). Now the joint posterior for $u$ and $\sigma^2$ given the data $m$ reads:
\begin{eqnarray}
\fl \quad p_{post}(u,\sigma^2|m) \propto \nonumber \\
\fl \qquad \left( \frac{1}{2 \sigma^2} \right)^{\case{k}{2}}  \exp \left( -\frac{1}{2 \, \sigma^2} \| m - A \, u \|^{2}_{2} -  \lambda | D \, u| \right) \exp \left( -(\alpha +1) \log(\sigma^2) - \frac{\beta}{\sigma^2} \right)  \nonumber \\
\fl \quad = \exp \left( -\frac{\case{1}{2} \| m - A \, u \|^{2}_{2} + \beta}{\sigma^2}  -  \lambda | D \, u|  -(\alpha +1 +k/2) \log(\sigma^2)  \right)\label{eq:JoinPost}
\end{eqnarray}
A comparison with \eref{eq:InvGamma} shows that the conjugacy property of the prior has the effect that the posterior of $\sigma^2$ conditioned on both $m$ and $u$ is, again, an inverse gamma distribution with shape and scale parameters $\tilde{\alpha}$ and $\tilde{\beta}$ given by :
\begin{equation}
 \tilde{\alpha} = \alpha + k/2; \qquad \tilde{\beta} = \case{1}{2} \| m - A \, u \|^{2}_{2} + \beta
\end{equation}
This allows us to perform Bayesian inference for the joint posterior using Gibbs sampling: For sampling along a component of $u$ we can use the fast samplers presented here\footnote{The practical implementation has to be slightly adopted in the sense that the varying $\sigma^2$ has to be removed from all expressions that are precomputed and added back to them at runtime.}, and for sampling over $\sigma^2$ conditioned on all other parameters we can use standard implementations of gamma samplers. Such a Gibbs sampler can, e.g., be used to infer a joint CM estimate $(u_{jCM},\sigma^2_{jCM})$. The information given by $\sigma^2_{jCM}$ can be used to evaluate or improve the measurement setup or to inform other Bayesian reconstructions. The estimate $u_{jCM}$ compared to an estimate assuming a single, constant $\sigma^2$ contains the \emph{marginalized} uncertainty about $\sigma^2$ and may yield a more robust estimate of $u$ in practical applications. In principle, it is possible to marginalize over $\sigma^2$ explicitly (the computation is similar to \cite{Ge06}):
\begin{eqnarray}
\fl \quad p_{post}(u|m) &= \int p_{post}(u,\sigma^2|m)  \rmd \sigma^2 \propto \left( \case{1}{2} \| m - A \, u \|^{2}_{2} + \beta \right)^{-(\alpha + k/2)} p_{pr}(u)\\
\fl \quad &\propto \left( 1 + \frac{t(u)^2}{\nu} \right)^{-\case{1}{2}(\nu + 1)}p_{pr}(u), \\
\fl \quad {\rm where} \quad &\nu = 2 \alpha + 1+ k; \quad t(u) = \left( \frac{\case{1}{2} \| m - A \, u \|^{2}_{2}}{\beta  (2 \alpha + 1+ k)^{1/2}} \right)^{1/2}.
\end{eqnarray}
This is a (one sided) Student's t-distribution for $t(u)$ with $\nu$ \emph{degrees of freedom}. However, working with this distribution directly is more difficult since it is not \emph{log-concave} (this problem gets worse if an individual variance $\sigma^2_i$ for each noise channel is assumed). Working with the full joint posterior instead can circumvent some of these problems.\\
The Gibbs samplers presented here have to be generalized to work with arbitrary $D$, e.g., to deal with anisotropic total variation priors in arbitrary dimensions and with arbitrary boundary conditions. In addition, an extension to \emph{block sparse} priors, which rely on mixed L2-L1 norms, would be advantageous to, e.g., address isotropic total variation priors.\\ 
Parallelization of MCMC sampling is easily implemented. In the most basic form, $N$ independent chains are generated, each on one CPU. However, the efficiency of this approach is strongly limited by the burn-in and mixing time \cite{Li08}. If all chains are initialized at the same state, parallelization is only efficient, if the chains become independent very fast. Our results suggest that parallelization of Gibbs samplers will be way more efficient than of MH samplers.\\
The Gibbs sampling algorithms developed by us are fast enough to tackle sampling for Bayesian inversion techniques in real applications, which will be an important topic of future work. In many applications like, e.g., limited angle CT, exploring the full range of Bayesian inversion by also incorporating sample based analysis was, up to now, rather regarded as a theoretical option, see, e.g., \cite{SiKoJaKaKoLaPiSo03,KoSiJaKaKoLaPiSo03}. \\
In addition, theoretical questions concerning sparse Bayesian inversion, like, e.g., the ones addressed in \cite{LaSi04,LaSaSi09} can be also be addressed numerically (cf. Section \ref{subsubsec:VisRes1D}).




\ack
We would like to thank the anonymous referees and Martin Burger for their helpful critiques and comments that significantly improved this manuscript.

\appendix

\section{Code} \label{sec:App_Code}

On the authors homepage\footnote{Currently: http://wwwmath.uni-muenster.de/num/burger/organization/lucka} Matlab code supporting this publication is provided. It contains scripts to create the scenarios examined in the numerical studies as well as implementations of all the Gibbs sampling algorithms presented here. One should, however, mention that Matlab is not very well suited for these implementation as the sampling algorithms consist of very sequential but rather basic procedures. We therefore also provide alternative implementations of the samplers using Fortran within \verb|.mex|-files, which can be compared to the corresponding \verb|.m|-files. The speed-up factor for the RnGibbs Sampler ranges from $43$ to $18$ using $n = 63$ or $n = 1023$ for the 1D scenario and is $1.35$ for the 2D scenario with $n = 261\,121$ (the computation of $b$ is by far the most expensive computational task, and little gain can be expected from a direct implementation compared to Matlab).

\section{Implementation}\label{sec:App_Imp}

In this section, we give details on how to implement formulas \eref{eq:cdf},  \eref{eq:leftinvcdf} and  \eref{eq:rightinvcdf}. The complementary error function and its inverse are difficult to handle numerically because there are no identities that allow to rescale or shift their evaluation to other intervals. For the applications we address, problems due to limited precision occur when formulas \eref{eq:cdf}, \eref{eq:leftinvcdf} and  \eref{eq:rightinvcdf} are implemented directly (formula \eref{eq:cdf} is only required for applying ordered overrelaxation). Dependent on the signs of $\alpha_+$ and  $\alpha_-$, we use different alternative formulas that allow for a stable numerical evaluation. Additionally, we express ${\rm erfc}(x)$ in terms of the scaled complementary error function ${\rm erfcx}(x) = \exp(x^2) {\rm erfc}(x)$, which decays less fast for $x \rightarrow +\infty$. We only list the results here (the corresponding transformations are elementary but lengthy to write down). Because $c \geqslant 0$, not both $\alpha_+$ and  $\alpha_-$  can be negative, which leaves three different cases to examine: \\
\underline{$\alpha_+ > 0$,  $\alpha_- > 0$:} Let $\gamma_{++} := \left[ {\rm erfcx} \left( \alpha_+ \right) + {\rm erfcx} \left( \alpha_- \right) \right] $. Then, the parts of \eref{eq:cdf} are given by:  
\begin{eqnarray}
\fl \qquad \qquad \; y < 0 &:  \exp \left( - a y^2 + 2 \sqrt{a} y \alpha_+  \right)  {\rm erfcx}(-\sqrt{a} y + \alpha_+ )/\gamma_{++} \label{eq:cdf++<} \\
\fl \qquad \qquad \; y > 0 &: 1 - \exp \left( - a y^2 - 2 \sqrt{a} y \alpha_-  \right)  {\rm erfcx}(\sqrt{a} y + \alpha_- )/\gamma_{++} \label{eq:cdf++>}
\end{eqnarray}
The arguments of ${\rm erfcinv}$ in \eref{eq:leftinvcdf} and \eref{eq:rightinvcdf} are given by:
\begin{eqnarray}
\fl \quad {\rm In} \; \eref{eq:leftinvcdf} &: \quad  r  \exp\left(  - \alpha^2_+ \right)    \gamma_{++}   \label{eq:Aleftinvcdf++}\\
\fl \quad {\rm In} \; \eref{eq:rightinvcdf} &: \quad (1-r) \exp(-\alpha_-^2)  \gamma_{++}  \label{eq:Arightinvcdf++}
\end{eqnarray}
\underline{$\alpha_+ < 0$,  $\alpha_- > 0$:} Since erfcx increases very fast for $x \rightarrow -\infty$ one has to use the identity ${\rm erfcx}(-x) = 2\, \exp(x^2) - {\rm erfcx}(x)$. Let $\gamma_{-+} := \left[  {\rm erfcx} \left( - \alpha_+ \right) - {\rm erfcx} \left( \alpha_- \right) \right]$. Then, the formulas to implement the parts of \eref{eq:cdf} are given by:  
\begin{eqnarray}
 \fl \quad \begin{array}{r}
      y < 0 \\
      - \sqrt{a} y + \alpha_+ > 0
     \end{array}
&: \quad \frac{\exp \left[ - \left(  \sqrt{a} y - \alpha_+ \right)^2  \right] {\rm erfcx}(-\sqrt{a} y + \alpha_+ ) }{2 - \exp\left( - \alpha_+^2\right) \gamma_{-+} }    \label{eq:cdf-+<>}\\ 
 \fl \quad \begin{array}{r}
      y < 0 \\
      - \sqrt{a} y + \alpha_+ < 0
     \end{array}
&: \quad \frac{  \left\lbrace 2 - \exp \left[ - \left(  \sqrt{a} y - \alpha_+ \right)^2  \right] {\rm erfcx}(\sqrt{a} y - \alpha_+ )  \right\rbrace}{2 - \exp\left( - \alpha_+^2\right) \gamma_{-+}} \label{eq:cdf-+<<}   \\
 \fl \;\;\; \quad \begin{array}{r}
      y > 0 \\
      \sqrt{a} y + \alpha_- > 0
     \end{array}
&: \quad \quad 1 - \frac{   \exp \left[ - \left(  \sqrt{a} y - \alpha_- \right)^2  \right] {\rm erfcx}(\sqrt{a} y + \alpha_- ) }{2 \exp \left( \frac{bc}{a}\right)  - \exp\left( - \alpha_-^2\right) \gamma_{-+}} \label{eq:cdf-+>>}\\
 \fl \;\;\; \quad \begin{array}{r}
      y > 0 \\
      \sqrt{a} y + \alpha_- < 0
     \end{array}
&: \quad 1 - \frac{2 -   \exp \left[ - \left(  \sqrt{a} y - \alpha_- \right)^2  \right] {\rm erfcx}(-\sqrt{a} y - \alpha_- ) }{2 \exp \left( \frac{bc}{a}\right)  - \exp\left( - \alpha_-^2\right) \gamma_{-+}} \label{eq:cdf-+><}
\end{eqnarray}
The arguments of ${\rm erfcinv}$ in \eref{eq:leftinvcdf} and \eref{eq:rightinvcdf} are given by:
\begin{eqnarray}
\fl \quad {\rm In} \; \eref{eq:leftinvcdf} &: \quad   r \left[  2 - \exp(-\alpha_+^2)  \gamma_{-+} \right]  \label{eq:Aleftinvcdf-+} \\
\fl \quad {\rm In} \; \eref{eq:rightinvcdf} &: \quad   (1-r) \left[  2 \exp \left( \frac{bc}{a} \right)   - \exp(-\alpha_-^2)  \gamma_{-+} \right]    \label{eq:Arightinvcdf-+} 
 \end{eqnarray}
\underline{$\alpha_+ > 0$,  $\alpha_- < 0$:} Let $\gamma_{+-} := \left[  {\rm erfcx} \left(  \alpha_+ \right) - {\rm erfcx} \left( - \alpha_- \right) \right]$. Then, the parts of \eref{eq:cdf} are given by:
\begin{eqnarray}
  \fl \qquad \qquad \; y < 0 &: \quad \frac{   \exp \left[ - \left(  \sqrt{a} y - \alpha_+ \right)^2  \right] {\rm erfcx}(-\sqrt{a} y + \alpha_+ ) }{2 \exp \left(- \frac{bc}{a}\right)  + \exp\left( - \alpha_+^2\right) \gamma_{+-}} \label{eq:cdf+-<}\\
 \fl \quad \begin{array}{r}
      y > 0 \\
      \sqrt{a} y + \alpha_- > 0
     \end{array}
&: \quad 1 - \frac{   \exp \left[ - \left(  \sqrt{a} y + \alpha_- \right)^2  \right] {\rm erfcx}(\sqrt{a} y + \alpha_- ) }{2  + \exp\left( - \alpha_-^2\right) \gamma_{+-}} \label{eq:cdf+->>}\\
 \fl \quad \begin{array}{r}
      y > 0 \\
      \sqrt{a} y + \alpha_- < 0
     \end{array}
&: \quad 1 - \frac{ 2 -   \exp \left[ - \left(  \sqrt{a} y + \alpha_- \right)^2  \right] {\rm erfcx}(-\sqrt{a} y - \alpha_- ) }{2  + \exp\left( - \alpha_-^2\right) \gamma_{+-}} \label{eq:cdf+-><}\\
\end{eqnarray}
The arguments of ${\rm erfcinv}$ in \eref{eq:leftinvcdf} and \eref{eq:rightinvcdf} are given by:
\begin{eqnarray}
\fl \quad {\rm In} \; \eref{eq:leftinvcdf} &: \quad  r \left[  2 \exp \left( \frac{-bc}{a} \right) + \exp(-\alpha_+^2)  \gamma_{+-} \right]  \label{eq:Aleftinvcdf+-} \\
\fl \quad {\rm In} \; \eref{eq:rightinvcdf} &: \quad (1-r) \left[  2  + \exp(-\alpha_-^2)  \gamma_{+-} \right]   \label{eq:Arightinvcdf+-}
 \end{eqnarray}
Using the above expressions directly can still lead to stability issues, because very large numbers are often multiplied with very small numbers. It is preferable to compute the logarithms of the expressions first. For this, let $x > 0$, $(x +y) > 0$ then:
\begin{equation}
 \log(x+y) = \log(x) + \log(1 + {\rm sign}(y) \exp(\log(|y|) - \log(x)) 
\end{equation}
Using this identity we can compute the logarithms of expressions \eref{eq:cdf++<} - \eref{eq:Arightinvcdf+-}. We note that ${\rm sign}\left[ \pm {\rm erfcx}(\cdot)\right] =  \pm 1$.
\begin{eqnarray}
\fl \quad \log \left[  \eref{eq:cdf++<} \right] = \left( - a y^2 + 2 \sqrt{a} y \alpha_+  \right) + \log \left[  {\rm erfcx}(-\sqrt{a} y + \alpha_+ ) \right]  - \log \left( \gamma_{++} \right) \\
\fl \quad \log \left[1-\eref{eq:cdf++>} \right] = \left( - a y^2 - 2 \sqrt{a} y \alpha_-  \right) + \log \left[  {\rm erfcx}(\sqrt{a} y + \alpha_- ) \right]  - \log \left( \gamma_{++} \right) \\
\fl \quad \log \left[  \eref{eq:Aleftinvcdf++} \right] = \log(r) - \alpha_+^2 + \log(\gamma_{++}) \\
\fl \quad \log \left[  \eref{eq:Arightinvcdf++} \right] = \log(1-r) - \alpha_-^2 + \log(\gamma_{++}) \\
\fl \quad \log \left[  \eref{eq:cdf-+<>} \right] = - \left( - \sqrt{a} y +  \alpha_+  \right)^2 + \log \left[  {\rm erfcx}(-\sqrt{a} y + \alpha_+ ) \right] - \log(2) \nonumber \\
\fl \qquad \quad \qquad - \log \left\lbrace  1 -  {\rm sign}(\gamma_{-+}) \exp \left[ -\alpha_+^2 + \log(|\gamma_{-+}|) - \log(2) \right] \right\rbrace \\
\fl \quad \log \left[  \eref{eq:cdf-+<>} \right] =  \log\left(  1 - \exp \left\lbrace  \log\left[ {\rm erfcx}(\sqrt{a} y - \alpha_+ ) \right] - \log\left( 2\right) - \left( \sqrt{a} y - \alpha_+ \right)^2  \right\rbrace    \right)  \nonumber \\
\fl \qquad \quad \qquad - \log \left\lbrace  1 -  {\rm sign}(\gamma_{-+}) \exp \left[ -\alpha_+^2 + \log(|\gamma_{-+}|) - \log(2) \right] \right\rbrace \\
\fl \quad \log \left[1-\eref{eq:cdf-+>>} \right] = - \left( \sqrt{a} y + 2 \alpha_- \right)^2 + \log \left[  {\rm erfcx}(\sqrt{a} y + \alpha_- ) \right] - \log(2) - \frac{bc}{a} \nonumber \\ 
\fl \qquad \quad \qquad - \log \left\lbrace  1 -  {\rm sign}(\gamma_{-+}) \exp \left[ -\alpha_-^2 + \log(|\gamma_{-+}|) - \log(2) - \frac{bc}{a} \right] \right\rbrace    \\
\fl \quad \log \left[1-\eref{eq:cdf-+><} \right] =  \log\left(  1 - \exp \left\lbrace  \log\left[ {\rm erfcx}(-\sqrt{a} y - \alpha_- ) \right] - \log\left( 2\right) - \left( \sqrt{a} y + \alpha_- \right)^2  \right\rbrace    \right)  \nonumber \\ 
\fl \qquad \quad  - \frac{bc}{a} - \log \left\lbrace  1 -  {\rm sign}(\gamma_{-+}) \exp \left[ -\alpha_-^2 + \log(|\gamma_{-+}|) - \log(2) - \frac{bc}{a} \right] \right\rbrace    \\
\fl \quad \log \left[  \eref{eq:Aleftinvcdf-+} \right]  = \log(r) +\log(2) \nonumber \\
\fl \qquad \quad \qquad + \log \left\lbrace  1 -  {\rm sign}(\gamma_{-+}) \exp \left[ -\alpha_+^2 + \log(|\gamma_{-+}|) - \log(2) \right] \right\rbrace    \\
\fl \quad \log \left[  \eref{eq:Arightinvcdf-+} \right]  = \log(1-r) +\log(2) + \frac{bc}{a} \nonumber \\ 
\fl \qquad \quad \qquad + \log \left\lbrace  1 -  {\rm sign}(\gamma_{-+}) \exp \left[ -\alpha_-^2 + \log(|\gamma_{-+}|) - \log(2) - \frac{bc}{a} \right] \right\rbrace    \\
\fl \quad \log \left[  \eref{eq:cdf+-<} \right] = - \left( - \sqrt{a} y +  \alpha_+  \right)^2 + \log \left[  {\rm erfcx}(-\sqrt{a} y + \alpha_+ ) \right] -\log(2) + \frac{bc}{a} \nonumber \\ 
\fl \qquad \quad \qquad  - \log \left\lbrace  1 + {\rm sign}(\gamma_{+-}) \exp \left[ -\alpha_+^2 + \log(|\gamma_{+-}|) - \log(2) + \frac{bc}{a} \right] \right\rbrace    \\
\fl \quad \log \left[  \eref{eq:cdf+->>} \right] = - \left(  \sqrt{a} y +  \alpha_-  \right)^2 + \log \left[  {\rm erfcx}(\sqrt{a} y + \alpha_- ) \right] -\log(2)  \nonumber \\
\fl \qquad \quad \qquad - \log \left\lbrace  1 + {\rm sign}(\gamma_{+-}) \exp \left[ -\alpha_-^2 + \log(|\gamma_{+-}|) - \log(2) \right] \right\rbrace\\
\fl \quad \log \left[  \eref{eq:cdf+-><} \right] = \log\left(  1 - \exp \left\lbrace  \log\left[ {\rm erfcx}(-\sqrt{a} y - \alpha_- ) \right] - \log\left( 2\right) - \left( \sqrt{a} y + \alpha_- \right)^2  \right\rbrace  \right) \nonumber \\
\fl \qquad \quad \qquad - \log \left\lbrace  1 + {\rm sign}(\gamma_{+-}) \exp \left[ -\alpha_-^2 + \log(|\gamma_{+-}|) - \log(2) \right] \right\rbrace\\
\fl \quad \log \left[  \eref{eq:Aleftinvcdf+-} \right]  = \log(r) +\log(2) - \frac{bc}{a} \nonumber \\ 
\fl \qquad \quad \qquad + \log \left\lbrace  1 + {\rm sign}(\gamma_{+-}) \exp \left[ -\alpha_+^2 + \log(|\gamma_{+-}|) - \log(2) + \frac{bc}{a} \right] \right\rbrace    \\
\fl \quad \log \left[  \eref{eq:Arightinvcdf+-} \right]  = \log(1-r) +\log(2)  \nonumber \\
\fl \qquad \quad \qquad + \log \left\lbrace  1 + {\rm sign}(\gamma_{+-}) \exp \left[ -\alpha_-^2 + \log(|\gamma_{+-}|) - \log(2) \right] \right\rbrace   
\end{eqnarray}
Now, for \eref{eq:leftinvcdf} and \eref{eq:rightinvcdf}, if $w$ denotes the logarithm of the argument of ${\rm erfcinv}$, one can compute ${\rm erfcinvlog}(w) :=  {\rm erfcinv}\left[  \exp ( w ) \right] $ using a standard implementation of erfcinv if $w$ is not too small (the loss of precision using $\exp(w)$ instead of computing the full argument of erfcinv is negligible since the variation of erfcinv is very small even on logarithmic scale). However, even using 64 bit precision is not sufficient for the applications we address.  Therefore, we use an asymptotic approximation of ${\rm erfcinvlog}(w)$ for $w < -680$ from \cite{dlmf11}.\\
An approximation of $z = {\rm erfcinv}\left[ \exp(w) \right] $ for $w \longrightarrow - \infty$ is given by:
\begin{eqnarray}
        \theta &:= - \log(\pi) - \log(-w) \nonumber \\
        v &:= (-\theta - 2) \nonumber \\
        s &:= 2/(\theta - 2  w) \nonumber \\
        a_2 &:= \frac{1}{8}  v \nonumber \\
        a_3 &:= - \frac{1}{32} (v^2 + 6v - 6) \nonumber \\
        a_4 &:= \frac{1}{384} (4 v^3 + 27 v^2 + 108 v - 300) \nonumber \\
        z  &\approx s^{-1/2} + a_2  s^{3/2} + a_3  s^{5/2} + a_4  s^{7/2} 
\end{eqnarray}
The discrepancy of this approximation to the implementation of erfcinv in Matlab is $2.34 \cdot 10^{-12}$ for $w = - 690$  and as it is an asymptotic formula, the error further decreases for $w \rightarrow - \infty$.


\section*{References}

\bibliographystyle{abbrv}
\bibliography{all}

\end{document}